\documentclass[reqno]{amsart}

\usepackage{amsmath,amsfonts,amssymb,amsthm}
\usepackage{mathtools}
\usepackage{mathrsfs}
\usepackage[inline]{enumitem}
\usepackage{todonotes}

\makeatletter
\newcommand{\inlineitem}[1][]{%
	\ifnum\enit@type=\tw@
	{\descriptionlabel{#1}}
	\hspace{\labelsep}%
	\else
	\ifnum\enit@type=\z@
	\refstepcounter{\@listctr}\fi
	\quad\@itemlabel\hspace{\labelsep}%
	\fi}
\makeatother

\makeatletter
\g@addto@macro{\endabstract}{\@setabstract}
\newcommand{\authorfootnotes}{\renewcommand\thefootnote{\@fnsymbol\c@footnote}}%
\makeatother

\theoremstyle{plain}
\newtheorem{theorem}{Theorem}[section]
\newtheorem{corollary}[theorem]{Corollary}
\newtheorem{proposition}[theorem]{Proposition}
\newtheorem{lemma}[theorem]{Lemma}
\theoremstyle{definition}
\newtheorem{definition}[theorem]{Definition}

\newtheorem{remark}[theorem]{Remark}
\newtheorem{example}[theorem]{Example}

\newcommand{\cG}{\mathcal G}
\newcommand{\cT}{\mathcal T}

\newcommand{\bd}{\mathbf d}
\newcommand{\be}{\mathbf e}
\newcommand{\bx}{\mathbf x}
\newcommand{\by}{\mathbf y}

\usepackage[export]{adjustbox}
\usepackage{forest}
\usetikzlibrary{shapes.geometric,decorations.markings}
\forestset{
	my tree/.style={
		before typesetting nodes={
			for tree={
				if={>O+tt=!O_=&{content}{}{n}{1}}{
					label/.process={Ow}{content}{left:##1},
				}{
					label/.process={Ow}{content}{right:##1},
				},
				content=,
			},
		},
		where level=0{
			draw,
			baseline,
			fill,
			circle,
		}{
			draw,
			fill,
			circle,
		},
		for tree={inner sep=1.3pt, l sep+=1pt, s sep'+=15pt},
	},
	default preamble={my tree},
	cut/.style={
		tikz+={
			\path ()  -- (!u);
		}
	}
}

\begin{document}
	
	\begin{center}
		\LARGE 
		Families of graphs with twin pendent paths and the Braess edge \par \bigskip
		
		\normalsize
		\authorfootnotes
		Sooyeong Kim\footnote{Department of Energy, Systems, Territory and Construction Engineering, Universit{\`a} di Pisa (sooyeong.kim@ing.unipi.it)} \par \bigskip
		
		\today
	\end{center}
\begin{abstract}
	In the context of a random walk on an undirected graph, Kemeny's constant can measure the average travel time for a random walk between two randomly chosen vertices. We are interested in graphs that behave counter-intuitively in regard to Kemeny's constant: in particular, we examine graphs with a cut-vertex at which at least two branches are paths, regarding whether the insertion of a particular edge into a graph results in an increase of Kemeny's constant. We provide several tools for identifying such an edge in a family of graphs and for analising  asymptotic behaviour of the family regarding the tendency to have that edge; and classes of particular graphs are given as examples. Furthermore, asymptotic behaviours of families of trees are described.
\end{abstract}



\section{Introduction}
%

Imagine a situation that adding roads to a road network in order to reduce traffic congestion results in, on contrary to one's expectation, slowing down overall traffic flow (called Braess's paradox \cite{Braess:paradox}). Random walks on graphs can also exhibit a version of this paradox. The family of random walks on undirected graphs is a special type of Markov chain: the transition probability from an initial state to another is given by the inverse of the degree of the vertex corresponding to the initial state. The parameter known as Kemeny's constant can be used to measure the average time for travel of a Markov chain between two randomly chosen states; so, in the context of random walks, it can be interpreted as a measure of how
well-connected the vertices of a graph are. Related applications can be found in \cite{Emma:Kemeny} for detecting potential super-spreaders of COVID-19, and in \cite{Crisostomi:Google} for determining `critical' roads in vehicle traffic networks based on Markov chains.

Kemeny's constant can serve as a proxy for identifying an edge exhibiting the version of the paradox \cite{KirklandZeng}, by examining an edge whose insertion into an undirected graph corresponding to a road network increases Kemeny's constant for random walks on the graph (such an edge is called a \textit{Braess edge}) that corresponds to travel times on the network. In the present paper, we study under what circumstances graphs can have a Braess edge in order to see what type of graphs exhibit the version of the paradox.

The term `Braess edge' is introduced in \cite{KirklandZeng}, and acknowledges Dietrich Braess who studies Braess's paradox for traffic networks \cite{Braess:paradox}. Kirkland and Zeng \cite{KirklandZeng} provides a particular family of trees, with a vertex adjacent to two pendent vertices (such two vertices are called \textit{twin pendent vertices}), such that inserting an edge between the twin pendent vertices causes Kemeny's constant to increase. Furthermore, Ciardo \cite{CiardoparadoxicalTwins} extends the result to all connected graphs with twin pendent vertices. Unlike the works \cite{KirklandZeng} and \cite{CiardoparadoxicalTwins}, Hu and Kirkland \cite{HuKirkland} establishes equivalent conditions for complete multipartite graphs and complete split graphs to have every non-edge as a Braess edge.

Our work is to generalise the circumstances in \cite{KirklandZeng,CiardoparadoxicalTwins} where graphs have a pair of twin pendent vertices; so, we consider graphs that can be constructed from a connected graph and two paths by identifying a vertex of the graph and a pendent vertex of each path. We call the two paths \textit{twin pendent paths} in the constructed graph. In Section \ref{Section:Braess edges on pendent}, a formula is derived that identifies a graph with twin pendent paths in which the non-edge between the pendent vertices of the twin pendent paths is a Braess edge. In Section \ref{Section: Asymptotic}, a combinatorial expression is provided in order to investigate asymptotic behaviour of a family of graphs with twin pendent paths regarding the tendency to have the non-edge, between the pendent vertices of the twin pendent paths, as a Braess edge. Furthermore, several families of graphs are discussed throughout Sections \ref{Section:Braess edges on pendent}, \ref{Section: Asymptotic}, and \ref{Section:Asymptotic for trees}. In particular, asymptotic behaviours of families of trees are characterized in Section \ref{Section:Asymptotic for trees}.

\section{Preliminaries}

Throughout the paper, we assume all graphs to be simple and undirected.

We shall introduce necessary terminology and notation in graph theory. Let $G$ be a graph of order $n$ with vertex set $V(G)$ and edge set $E(G)$ where $n=|V(G)|$. An edge joining vertices $v$ and $w$ of $G$ is denoted by $v\sim w$. Let $m_G$ be defined as $|E(G)|$. The subgraph of $G$ \textit{induced} by a subset $S$ of $V(G)$ is the graph with vertex set $S$, where two vertices in $S$ are adjacent if and only if they are adjacent in $G$. For $v\in V(G)$, we denote by $\mathrm{deg}_G(v)$ the degree of $v$. A vertex $v$ of a graph $G$ is said to be \textit{pendent} if $\mathrm{deg}_G(v)=1$. Given a labelling of $V(G)$, we define $\mathbf{d}_G$ to be the column vector whose $i^\text{th}$ component is $\mathrm{deg}_G(v_i)$ for $1\leq i\leq n$, where $v_i$ is the $i^\text{th}$ vertex in $V(G)$. For $v,w\in V(G)$, the distance between $v$ and $w$ in $G$ is denoted by $\mathrm{dist}_G(v,w)$. For a connected graph $G$ with a vertex $v$, the \textit{eccentricity} $e_G(v)$ of $v$ in $G$ is $e_G(v)=\mathrm{max}\{\mathrm{dist}_G(v,w)|w\in V(G)\}$. The \textit{diameter}, denoted $\mathrm{diam}(G)$, of $G$ is $\mathrm{diam}(G)=\mathrm{max}\{e_G(v)|v\in V(G)\}$.

The \textit{trivial} graph is the graph of order $1$. A \textit{tree} is a connected graph that has no cycles. A \textit{forest} is a graph whose connected components are trees. A \textit{spanning tree} (resp. a \textit{spanning forest}) of $G$ is a subgraph that is a tree (resp. a forest) and includes all of the vertices of $G$. A \textit{$k$-tree spanning forest} of $G$ is a spanning forest that consists of $k$ trees. For $v\in V(G)$, we use $G-v$ to denote the graph obtained from $G$ by the deletion of $v$. A vertex $v$ of a connected graph $G$ is called a \textit{cut-vertex} of $G$ if $G-v$ is disconnected. If $G-v$ has $k$ connected components $G_1,\dots,G_k$ for some $k\geq 2$, then the subgraph induced by $V(G_i)\cup\{v\}$ for $1\leq i\leq k$ is called a \textit{branch} of $G$ at $v$.

Let us introduce several types of connected graphs. We denote the complete graph of order $n$ by $K_n$, the cycle of length $n$ by $C_n$, and the path on $n$ vertices by $P_n$. If we need to specify the ordering of vertices of a cycle or a path, then we use $C_n=(v_1,v_2,\dots,v_n,v_1)$ to denote the cycle whose vertices are labelled by $v_1,\dots,v_n$, and whose edges are $v_1\sim v_n$ and $v_i\sim v_{i+1}$ for $i=1,\dots,n-1$; similarly, $P_n=(v_1,v_2,\dots,v_n)$ denotes the path whose vertices are labelled by $v_1,\dots,v_n$ and whose edges are $v_i\sim v_{i+1}$ for $i=1,\dots,n-1$. A \textit{star} $S_n$ is a tree on $n$ vertices with one vertex of degree $n-1$. For $n\geq 3$, $v$ is called the \textit{centre} vertex of $S_n$ if $\mathrm{deg}_{S_n}(v)=n-1$. For $n>k\geq 1$, a \textit{broom} $\mathcal{B}_{n,k}$ is a tree constructed from the path on $k$ vertices by adding $n-k$ pendent vertices to one pendent vertex on the path.

Let $k\geq 2$, and let $G_i$ be a graph with $v_i\in V(G_i)$ for $i=1,\dots,k$. Suppose that $V(G_1),\dots,V(G_k)$ are disjoint. Let $v\notin V(G_i)$ for $i=1,\dots,k$. We consider a graph $G$ with vertex set $V(G)=\{v\}\cup\left(\bigcup_{i=1}^k \left(V(G_i)-v_i\right)\right)$, where two vertices $x$ and $y$ in $G$ are adjacent if and only if it satisfies one of the following: \begin{enumerate*}[label=(\roman*)]
	\item $x\sim y\in \bigcup_{i=1}^kE(G_i)$; and 
	\item one of $x$ and $y$ is $v$, and the other is a vertex adjacent to $v_j$ in $G_j$ for some $1\leq j\leq k$.
\end{enumerate*}
Then, we say that the graph $G$ is obtained from $G_1,\dots, G_k$ by \textit{identifying} vertices $v_1,\dots,v_k$ as $v$. 

Let $G$ be a graph. Let $P_{k_1+1}=(v_1,\dots,v_{k_1+1})$ and $P_{k_2+1}=(w_1,\dots,w_{k_2+1})$ where $k_1$ and $k_2$ are non-negative integers with $k_1+k_2\geq 2$. Suppose that $\widetilde{G}$ is the graph obtained from $G$, $P_{k_1+1}$, and $P_{k_2+1}$ by identifying a vertex $v$ of $G$, $v_1$, and $w_1$ as $v$. We say that the paths $(v,v_2,\dots,v_{k_1+1})$ and $(v,w_2\dots,w_{k_2+1})$ in $\widetilde{G}$ are \textit{twin pendent paths}. Then, the pendent vertices of the twin pendent paths in $\widetilde{G}$ are $v_{k_1+1}$ and $w_{k_2+1}$. We remark that considering the construction of $\widetilde{G}$, it is reasonable to assume that $k_1$ and $k_2$ are permitted to be zero (we consider it throughout this paper), as opposed to one's anticipation from the word `twin pendent paths', in that $P_{k_1+1}$ and $P_{k_2+1}$ both are of length at least $1$.

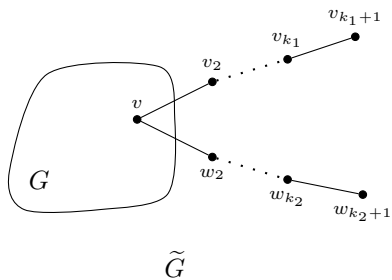
\begin{figure}[t!]
	\begin{center}
		\begin{tikzpicture}
		\tikzset{enclosed/.style={draw, circle, inner sep=0pt, minimum size=.10cm, fill=black}}
		\draw plot [smooth cycle] coordinates {(0.3,.3)(1.3,.3)(2.2,.5) (2.3,1.5)(2.1,2.2)(0.6,2.1)(0.1,0.8)} node at (0.5,0.7) {$G$};
		\node[label={below, yshift=0cm: $\widetilde{G}$}] (tilG) at (2.3,0) {};
		\begin{scriptsize}
		\node[enclosed, label={above, yshift=0cm: $v$}] (v) at (1.8,1.5) {};
		\node[enclosed, label={above, yshift=0cm: $v_2$}] (v1) at (2.8,2) {};
		\node[enclosed, label={above, yshift=0cm: $v_{k_1}$}]  (v2) at (3.8,2.3) {};
		\node[enclosed, label={above, yshift=0cm: $v_{k_1+1}$}] (v3) at (4.7,2.6) {};
		\node[enclosed, label={below, yshift=0cm: $w_2$}] (w1) at (2.8,1) {};
		\node[enclosed, label={below, yshift=0cm: $w_{k_2}$}]  (w2) at (3.8,0.7) {};
		\node[enclosed, label={below, yshift=0cm: $w_{k_2+1}$}] (w3) at (4.8,0.5) {};
		\end{scriptsize}
		\draw (v) -- (v1);
		\draw[thick, loosely dotted] (v1) -- (v2);
		\draw (v2) -- (v3);
		\draw (v) -- (w1);
		\draw[thick, loosely dotted] (w1) -- (w2);
		\draw (w2) -- (w3);
		\end{tikzpicture}
	\end{center}
	\caption{An illustration of twin pendent paths in $\widetilde{G}$.}\label{Figure:twinpendentpaths}
\end{figure}

Consider a discrete, finite, time-homogeneous Markov chain whose finite state space is $\{1,\dots,n\}$. The Markov chain can be represented by the transition matrix $M$. (We refer the interested reader to \cite{Seneta:Markov} for the necessary background for Markov chains.) Then, Kemeny's constant $\kappa(M)$ is defined as $\sum_{j\neq i}^n m_{i,j}w_j$, where $m_{i,j}$ is the mean first passage time from state $i$ to state $j$, and $w_j$ is the $j^\text{th}$ entry of the stationary distribution. Note that Kemeny's constant is independent of $i$. It is found in \cite{Levene:KemenyInterpre} that $\kappa(M)+1=\sum_{i=1}^n\sum_{j=1}^n w_im_{i,j}w_j$. This admits the interpretation of Kemeny's constant in terms of the expected number of steps from a randomly-chosen initial state to a randomly-chosen final state. Alternatively, $\kappa(M)$ can be expressed as $\kappa(M)=\sum_{j=2}^{n}\frac{1}{1-\lambda_j}$ where $1,\lambda_2,\dots,\lambda_n$ are the eigenvalues of $M$. For the details, the reader may refer to \cite{KemenySnell}.

For our work, we use the combinatorial expression for Kemeny's constant for a random walk on a connected graph in \cite{KirklandZeng}. In order to emphasize that we are dealing with random walks on connected graphs, given a connected graph $G$, we use $\kappa(G)$ to denote Kemeny's constant for the transition matrix of the random walk on $G$. We denote by $\tau_{G}$ the number of spanning trees of $G$, and by $\mathcal{F}_{G}(i;j)$ the set of $2$-tree spanning forests of $G$ such that one of the two trees contains a vertex $i$ of $G$, and the other has a vertex $j$ of $G$. We define $F_{G}$ to be the matrix given by $F_{G}=[f_{i,j}^{G}]$ where $f_{i,j}^{G}=|\mathcal{F}_{G}(i;j)|$. Note that $f_{i,i}^{G}=0$, that is, the diagonal entries of $F_G$ are zero. Recall that $m_G$ is the number of edges of $G$. Then, Kemeny's constant for the transition matrix of the random walk on $G$ is given by
\begin{equation*}
\kappa(G)=\frac{\mathbf{d}_G^T F_G\mathbf{d}_G}{4m_G\tau_G}.
\end{equation*}
A non-edge $e$ of $G$ is said to be a \textit{Braess edge} for $G$ if $\kappa(G)<\kappa(G\cup e)$ where $G\cup e$ is the graph obtained from $G$ by adding $e$ to $G$. A connected graph $G$ is said to be \textit{paradoxical} \cite{CiardoparadoxicalTwins} if there exists a Braess edge for $G$.

We also introduce some useful notation. We denote by $\mathbf{1}_k$ the all ones vector of size $k$, by $\mathbf{0}_k$ the all zeros vector of size $k$, and by $J_{p,q}$ the all ones matrix of size $p\times q$. If $k=p=q$, we write $J_{p,q}$ as $J_k$. The subscripts $k$ and a pair of $p$ and $q$ are omitted if their sizes are clear from the context. We denote by $\be_v$ the column vector whose component in $v^\text{th}$ position is $1$ and zeros elsewhere, and denote by $\mathbf{f}_G^v$ the $v^{\text{th}}$ column of $F_{G}$. Then, the $i^{\text{th}}$ entry of $\mathbf{f}_G^v$ is $f_{i,v}^{G}$. Moreover, $\mathbf{f}_G^v$ can be written as $\mathbf{f}^v$ if $G$ is clear from the context.

We assume familiarity with basic material on graph theory. We refer the reader to \cite{book:GraphAndDigraphs} for the necessary background. In what follows, we omit $G$ that is a subscript or a superscript in the notation described in this section when no confusion arises, and we use boldface lowercase letters to denote column vectors.

\section{Graphs with twin pendent paths and the Braess edge}\label{Section:Braess edges on pendent}

In this section, for a connected graph with twin pendent paths, we provide an equivalent condition for the non-edge between the pendent vertices of the twin pendent paths to be a Braess edge. Moreover, we examine several families of graphs through the equivalent condition.

We begin with investigating the components in the expression for Kemeny's constant $\kappa(G)$ where $G$ is a connected graph with a cut-vertex---a connected graph with twin pendent paths has a cut-vertex.

\begin{proposition}\label{gen:H1andH2}
	Let $H_1$ and $H_2$ be connected graphs, and let $v_1\in V(H_1)$ and $v_2\in V(H_2)$. Assume that $G$ is obtained from $H_1$ and $H_2$ by identifying $v_1$ and $v_2$ as a vertex $v$. Suppose that $\widetilde{H}_1=H_1-v_1$ and $\widetilde{H}_2=H_2-v_2$. Then, labelling the vertices of $G$ in order of $V(\widetilde{H}_1)$, $v$, and $V(\widetilde{H}_2)$, we have:
	\begin{align}\nonumber
	\mathbf{d}_{G}^T&=[\mathbf{d}_{H_1}^T\;\mathbf{0}^T_{|V(\widetilde{H}_2)|}]+[\mathbf{0}^T_{|V(\widetilde{H}_1)|}\;\mathbf{d}_{H_2}^T],\\\nonumber
	m_{G}&=m_{H_1}+m_{H_2},\\\nonumber
	\tau_{G}&=\tau_{H_1}\tau_{H_2},\\\label{F_G}
	F_{G}&=\left[\begin{array}{c|c|c}
	\tau_{H_2}F_{\widetilde{H}_1} & \tau_{H_2}\mathbf{f}_1 & \tau_{H_2}\mathbf{f}_1\mathbf{1}^T+\tau_{H_1}\mathbf{1}\mathbf{f}_2^T \\\hline
	\tau_{H_2}\mathbf{f}^T_1 & 0 & \tau_{H_1}\mathbf{f}^T_2\\\hline
	\tau_{H_1}\mathbf{f}_2\mathbf{1}^T+\tau_{H_2}\mathbf{1}\mathbf{f}_2^T & \tau_{H_1}\mathbf{f}_2  & \tau_{H_1}F_{\widetilde{H}_2}
	\end{array}\right],
	\end{align}
	where $\mathbf{f}_1$ and $\mathbf{f}_2$ are the column vectors obtained from $\mathbf{f}_{H_1}^v$ and $\mathbf{f}_{H_2}^v$ by deleting the $v^\text{th}$ component (which is $0$), respectively. Furthermore, we obtain 
	\begin{align}\label{dfd:expression}
	\mathbf{d}_G^TF_G\mathbf{d}_G
	=\tau_{H_2}\mathbf{d}_{H_1}^TF_{H_1}\mathbf{d}_{H_1}+\tau_{H_1}\mathbf{d}_{H_2}^TF_{H_2}\mathbf{d}_{H_2}+4\tau_{H_2}m_{H_2}\mathbf{d}_{H_1}^T\mathbf{f}_{H_1}^v+4\tau_{H_1}m_{H_1}\mathbf{d}_{H_2}^T\mathbf{f}_{H_2}^v.
	\end{align} 
\end{proposition}
\begin{proof}
	The conclusions for $\mathbf{d}_G$ and $m_G$ are readily established. We shall consider $F_G$ and $\tau_G$. Since $f_{i,j}^G=f_{j,i}^G$ for all $i,j\in V(G)$, $F_G$ is symmetric. Hence, we only need to verify the entries above the main diagonal. Note that $v$ is a cut-vertex of $G$. Since all spanning trees of $G$ can be obtained from spanning trees of $H_1$ and of $H_2$ by identifying $v_1$ and $v_2$ as $v$, we have $\tau_G=\tau_{H_1}\tau_{H_2}$. Let $i,j\in V(H_1)$. For each spanning forest of $H_1$ in $\mathcal{F}_{H_1}(i;j)$, we can obtain $\tau_{H_2}$ spanning forests of $G$ in $\mathcal{F}_{G}(i;j)$ from the forest of $H_1$ and each of $\tau_{H_2}$ spanning trees of $H_2$ by identifying $v_1$ and $v_2$. Therefore, $f_{i,j}^{G}=\tau_{H_2}f_{i,j}^{H_1}$ for $i,j\in V(H_1)$. Similarly, for $i,j\in V(H_2)$, we have $f_{i,j}^{G}=\tau_{H_1}f_{i,j}^{H_2}$. Let $i\in V(\widetilde{H}_1)$ and $j\in V(\widetilde{H}_2)$. The set $\mathcal{F}_{G}(i;j)$ is a disjoint union of $A_i$ and $A_j$, where $A_i$ is the set of spanning forests of $G$ in $\mathcal{F}_{G}(i;j)$ such that the tree having the vertex $i$ among the two trees contains $v$, and $A_j=\mathcal{F}_G(i;j)\backslash A_i$. Since for each spanning forest in $A_i$ the tree with $i$ has $v$, the tree contains a spanning tree of $H_1$ as a subtree. So, any forest in $A_i$ can be constructed from a spanning tree of $H_1$ and a spanning forest in $\mathcal{F}_{H_2}(v_2;j)$ with $v_1$ and $v_2$ identified as $v$. Hence, we have $|A_i|=\tau_{H_1}f_{v_2,j}^{H_2}$. Note that $f_{i,v_1}^{H_1}=f_{v_1,i}^{H_1}$. Applying an analogous argument to the case $|A_j|$, we have $|A_j|=\tau_{H_2}f_{i,v_1}^{H_1}$. Thus, $f_{i,j}^{G}=\tau_{H_2}f_{i,v}^{H_1}+\tau_{H_1}f_{v,j}^{H_2}$ for $i\in V(\widetilde{H}_1)$ and $j\in V(\widetilde{H}_2)$. Therefore, our desired results for $F_G$ and $\tau_G$ are obtained.
	
	Now, we shall prove \eqref{dfd:expression}. Labelling the vertices of $H_1$ (resp. $H_2$) in order of $V(\widetilde{H}_1)$ and $v$ (resp. $v$ and $V(\widetilde{H}_2)$), we have
	$$
	\tau_{H_2}F_{H_1}=\left[\begin{array}{cc}
	\tau_{H_2}F_{\widetilde{H}_1} & \tau_{H_2}\mathbf{f}_1  \\
	\tau_{H_2}\mathbf{f}^T_1 & 0 
	\end{array}\right],\;\tau_{H_1}F_{H_2}=\left[\begin{array}{cc}
	0 & \tau_{H_1}\mathbf{f}^T_2\\
	\tau_{H_1}\mathbf{f}_2  & \tau_{H_1}F_{\widetilde{H}_2}
	\end{array}\right].
	$$
	Note that $(\mathbf{f}_{H_1}^v)^T=\begin{bmatrix}
	\mathbf{f}_1^T & 0
	\end{bmatrix}$ and $(\mathbf{f}_{H_2}^v)^T=\begin{bmatrix}
	0 & \mathbf{f}_2^T
	\end{bmatrix}$. Then, 
	$$
	\begin{bmatrix}
	\tau_{H_2}\mathbf{f}_1 & \tau_{H_2}\mathbf{f}_1\mathbf{1}_{|V(\widetilde{H}_2)|}^T+\tau_{H_1}\mathbf{1}_{|V(\widetilde{H}_1)|}\mathbf{f}_2^T\\
	0 & \tau_{H_1}\mathbf{f}^T_2
	\end{bmatrix}=\tau_{H_2}\mathbf{f}_{H_1}^v\mathbf{1}_{|V(H_2)|}^T+\tau_{H_1}\mathbf{1}_{|V(H_1)|}(\mathbf{f}_{H_2}^v)^T.
	$$
	Considering $\mathbf{d}_{G}^T=[\mathbf{d}_{H_1}^T\;\mathbf{0}^T_{|V(\widetilde{H}_2)|}]+[\mathbf{0}^T_{|V(\widetilde{H}_1)|}\;\mathbf{d}_{H_2}^T]$ with $F_G$ in \eqref{F_G}, we have
	\begin{align*}\nonumber
	&\mathbf{d}_G^TF_G\mathbf{d}_G\\\nonumber=&\tau_{H_2}\mathbf{d}^T_{H_1}F_{H_1}\mathbf{d}_{H_1}+\tau_{H_1}\mathbf{d}^T_{H_2}F_{H_2}\mathbf{d}_{H_2}+2\mathbf{d}^T_{H_1}\begin{bmatrix}
	\tau_{H_2}\mathbf{f}_1 & \tau_{H_2}\mathbf{f}_1\mathbf{1}^T+\tau_{H_1}\mathbf{1}\mathbf{f}_2^T\\
	0 & \tau_{H_1}\mathbf{f}^T_2
	\end{bmatrix}\mathbf{d}_{H_2}
	\\\nonumber=&\tau_{H_2}\mathbf{d}_{H_1}^TF_{H_1}\mathbf{d}_{H_1}+\tau_{H_1}\mathbf{d}_{H_2}^TF_{H_2}\mathbf{d}_{H_2}+2\mathbf{d}_{H_1}^T\left(\tau_{H_2}\mathbf{f}_{H_1}^v\mathbf{1}_{|V(H_2)|}^T+\tau_{H_1}\mathbf{1}_{|V(H_1)|}(\mathbf{f}_{H_2}^v)^T\right)\mathbf{d}_{H_2}
	\\
	=&\tau_{H_2}\mathbf{d}_{H_1}^TF_{H_1}\mathbf{d}_{H_1}+\tau_{H_1}\mathbf{d}_{H_2}^TF_{H_2}\mathbf{d}_{H_2}+4\tau_{H_2}m_{H_2}\mathbf{d}_{H_1}^T\mathbf{f}_{H_1}^v+4\tau_{H_1}m_{H_1}\mathbf{d}_{H_2}^T\mathbf{f}_{H_2}^v.
	\end{align*}
\end{proof}

We can see from Proposition \ref{gen:H1andH2} that given $m_{H_i}$ and $\tau_{H_i}$ for $i=1,2$, $\mathbf{d}_G^TF_G\mathbf{d}_G$ can be computed from $\mathbf{d}_{H_i}^TF_{H_i}\mathbf{d}_{H_i}$ and $\mathbf{d}_{H_i}^T\mathbf{f}_{H_i}^v$ for $i=1,2$. The following examples regarding $K_n$, $C_n$, $P_n$ and $S_n$ present the corresponding quantities $\mathbf{d}^TF\mathbf{d}$ and $\mathbf{d}^T\mathbf{f}^v$, and assist us later to obtain several results and related examples.

\begin{example}\label{ex:K_n}
	{\rm Consider a complete graph $K_n$ on $n$ vertices. Then, $m=\binom{n}{2}$ and $\tau=n^{n-2}$ by Cayley's formula (see \cite{book:GraphAndDigraphs}). Note that $K_n$ is symmetric (see \cite{Godsil:AlgebraicGraph}), \textit{i.e.,} for any pair of edges of $K_n$, there is an automorphism that maps one edge to the other. So, $F_{K_n}=\alpha(J-I)$ where $\alpha=f_{i,j}^{K_n}$ for all $i,j\in V(K_n)$. Then, $\alpha$ is the determinant of the submatrix obtained from the Laplacian matrix of $K_n$ by deleting $i^\text{th}$ and $j^\text{th}$ rows and columns where $i\neq j$ (see \cite{Seth:MatrixTreeTheorem}). It can be seen that $\alpha=2n^{n-3}$. Therefore, for any vertex $v$ in $K_n$, we have}
	\begin{align*}
	&\mathbf{d}^TF\mathbf{d}=\alpha(n-1)^2\mathbf{1}^T(J-I)\mathbf{1}=2n^{n-2}(n-1)^3,\\
	&\mathbf{d}^T\mathbf{f}^v=\alpha(n-1)\mathbf{1}^T(\mathbf{1}-\be_v)=2n^{n-3}(n-1)^2.
	\end{align*}
\end{example}

\begin{example}\label{ex:C_n}
	{\rm Consider the cycle $C_n=(1,2,\dots,n,1)$ where $n\geq 3$. For $1\leq v \leq n$, we obtain
		\begin{align*}
		&F_{C_n}=\begin{bmatrix}
		\mathrm{dist}(i,j)(n-\mathrm{dist}(i,j))
		\end{bmatrix}_{1\leq i,j\leq n},\; \mathbf{d}=2\mathbf{1},\\
		&(\mathbf{f}^v)^T=\begin{bmatrix}
		(v-1)(n-(v-1)) & \cdots & 1\cdot(n-1) & 0 & 1\cdot(n-1) & \cdots & (n-v)v
		\end{bmatrix}.
		\end{align*}
		It can be checked that for $v=1,\dots,n$,
		\begin{align*}
		\mathbf{d}^TF\mathbf{d}=\frac{2}{3}(n-1)n^2(n+1)\;\text{and}\;\mathbf{d}^T\mathbf{f}^v=\frac{1}{3}(n-1)n(n+1).
		\end{align*}}
\end{example}

Note that for any tree $\cT$, $F_{\cT}$ is the distance matrix of $\cT$ (see \cite{KirklandZeng}), which is the matrix whose $(i,j)$-entry is the distance between $i$ and $j$.

\begin{example}\label{ex:P_n}
	{\rm Consider the path $P_n=(1,2,\dots,n)$ where $n\geq 2$. Let $v$ be a vertex of $P_n$. For $1\leq v \leq n$, we have \begin{align*}
		&F_{P_n}=\begin{bmatrix}
		|i-j|
		\end{bmatrix}_{1\leq i,j\leq n},\; \mathbf{d}=2\mathbf{1}_{n}-\mathbf{e}_1-\mathbf{e}_{n},\\
		&(\mathbf{f}^v)^T=\begin{bmatrix}
		v-1 & \cdots & 1 & 0 & 1 & \cdots & n-v
		\end{bmatrix}.
		\end{align*} 
		One can verify that for $v=1,\dots,n$,
		\begin{align*}
		&\mathbf{d}^TF_{}\mathbf{d}=4\mathbf{1}^TF\mathbf{1}-4\mathbf{1}^TF\be_1-4\mathbf{1}^TF\be_n+2\be_1^TF\be_n=\frac{4}{3}(n-1)^3+\frac{2}{3}(n-1),\\
		&\mathbf{d}^T\mathbf{f}^v=(v-1)^2+(n-v)^2.
		\end{align*}}
\end{example}

\begin{example}
	{\rm Consider a star $S_n$ of order $n$ where $n\geq 3$. Suppose that $n$ is the centre vertex. Then, we have
		\begin{align*}
		\mathbf{d}^T=\begin{bmatrix}
		\mathbf{1}_{n-1}^T & 0
		\end{bmatrix}+(n-1)\be_n,\;F_{S_n}=\begin{bmatrix}
		2(J-I) & \mathbf{1}_{n-1}\\
		\mathbf{1}_{n-1}^T & 0
		\end{bmatrix}.
		\end{align*}
		Hence, we have that for $v=1,\dots,n$,
		\begin{align*}
		&\mathbf{d}^TF\mathbf{d}=2\mathbf{1}_{n-1}^T(J-I)\mathbf{1}_{n-1}+2(n-1)^2=2(n-1)(2n-3),\\
		&\mathbf{d}^T\mathbf{f}^v=\begin{cases*}
		n-1, & \text{if $v=n$}\\
		3n-5, & \text{if $v\neq n$.}
		\end{cases*}
		\end{align*}}
\end{example}

We consider the following definition for clarity of exposition regarding our work in this paper.

\begin{definition}\label{Def:paradoxical}
	{\rm Let $G$ be a connected graph on $n$ vertices, and $v\in V(G)$. Fix two non-negative integers $k_1,k_2$ with $k_1+k_2\geq 2$. Let $\widetilde G(v,k_1,k_2)$ denote the graph obtained from $G$, $P_{k_1+1}=(v_1,\dots,v_{k_1+1})$ and $P_{k_2+1}=(w_1,\dots, w_{k_2+1})$ by identifying the vertices $v$, $v_1$ and $w_1$. Also, we denote by $\widehat G(v,k_1,k_2)$ the graph obtained from $\widetilde G(v,k_1,k_2)$ by inserting the edge $v_{k_1+1}\sim w_{k_2+1}$. We say that $G$ is {\em $(v,k_1,k_2)$-paradoxical} if $\kappa(\widehat G(v,k_1,k_2))>\kappa(\widetilde G(v,k_1,k_2))$. If $G$ is $(v,k_1,k_2)$-paradoxical for every $v\in V(G)$, then we say that $G$ is {\em $(k_1,k_2)$-paradoxical}.}
\end{definition}

Our main goal stated in the beginning of this section can be rephrased in terms of Definition \ref{Def:paradoxical}: given a connected graph $G$ with a vertex $v$, we shall find an equivalent condition for $G$ to be $(v,k_1,k_2)$-paradoxical.

\begin{example}\label{Ex:paradoxical}
	{\rm Consider the following graphs:
		
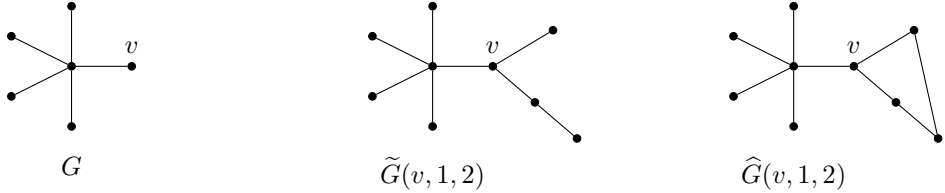
\begin{center}
			\begin{tikzpicture}[scale=0.80]
			\tikzset{enclosed/.style={draw, circle, inner sep=0pt, minimum size=.10cm, fill=black}}
			
			\node[enclosed, label={left, yshift=0cm: }] (v_1) at (-1,0.5) {};
			\node[enclosed, label={left, yshift=0cm: }] (v_2) at (0,-1) {};
			\node[enclosed, label={below, yshift=0cm:}] (v_3) at (0,0) {};
			\node[enclosed, label={right, yshift=0cm: }] (v_4) at (0,1) {};
			\node[enclosed, label={right, yshift=0cm: }] (v_5) at (-1,-0.5) {};
			\node[enclosed, label={above, yshift=0cm: $v$}] (v_6) at (1,0) {};
			\node[label={below, yshift=0cm: $G$}] (G) at (0,-1.2) {};
			\draw (v_2) -- (v_3);
			\draw (v_1) -- (v_3);
			\draw (v_3) -- (v_4);
			\draw (v_3) -- (v_5);
			\draw (v_3) -- (v_6);
			
			\begin{scope}[xshift=6cm]
			\tikzset{enclosed/.style={draw, circle, inner sep=0pt, minimum size=.10cm, fill=black}}
			
			\node[enclosed, label={left, yshift=0cm: }] (v_1) at (-1,0.5) {};
			\node[enclosed, label={left, yshift=0cm: }] (v_2) at (0,-1) {};
			\node[enclosed, label={below, yshift=0cm:}] (v_3) at (0,0) {};
			\node[enclosed, label={right, yshift=0cm: }] (v_4) at (0,1) {};
			\node[enclosed, label={right, yshift=0cm: }] (v_5) at (-1,-0.5) {};
			\node[enclosed, label={above, yshift=0cm: $v$}] (v_6) at (1,0) {};
			\node[enclosed, label={right, yshift=0cm: }] (v_9) at (2,0.6) {};
			\node[enclosed, label={right, yshift=0cm: }] (v_10) at (1.7,-0.6) {};
			\node[enclosed, label={right, yshift=0cm: }] (v_11) at (2.4,-1.2) {};
			\node[label={below, yshift=0cm: $\widetilde{G}(v,1,2)$}] (G) at (0,-1.2) {};
			\draw (v_2) -- (v_3);
			\draw (v_1) -- (v_3);
			\draw (v_3) -- (v_4);
			\draw (v_3) -- (v_5);
			\draw (v_3) -- (v_6);
			\draw (v_6) -- (v_9);
			\draw (v_6) -- (v_10);
			\draw (v_10) -- (v_11);
			\end{scope}
			
			\begin{scope}[xshift=12cm]
			\tikzset{enclosed/.style={draw, circle, inner sep=0pt, minimum size=.10cm, fill=black}}
			
			\node[enclosed, label={left, yshift=0cm: }] (v_1) at (-1,0.5) {};
			\node[enclosed, label={left, yshift=0cm: }] (v_2) at (0,-1) {};
			\node[enclosed, label={below, yshift=0cm:}] (v_3) at (0,0) {};
			\node[enclosed, label={right, yshift=0cm: }] (v_4) at (0,1) {};
			\node[enclosed, label={right, yshift=0cm: }] (v_5) at (-1,-0.5) {};
			\node[enclosed, label={above, yshift=0cm: $v$}] (v_6) at (1,0) {};
			\node[enclosed, label={right, yshift=0cm: }] (v_9) at (2,0.6) {};
			\node[enclosed, label={right, yshift=0cm: }] (v_10) at (1.7,-0.6) {};
			\node[enclosed, label={right, yshift=0cm: }] (v_11) at (2.4,-1.2) {};
			\node[label={below, yshift=0cm: $\widehat{G}(v,1,2)$}] (G) at (0,-1.2) {};
			\draw (v_2) -- (v_3);
			\draw (v_1) -- (v_3);
			\draw (v_3) -- (v_4);
			\draw (v_3) -- (v_5);
			\draw (v_3) -- (v_6);
			\draw (v_6) -- (v_9);
			\draw (v_6) -- (v_10);
			\draw (v_10) -- (v_11);
			\draw (v_9) -- (v_11);
			\end{scope}
			\end{tikzpicture}
		\end{center}
		With the aid of MATLAB\textsuperscript{\textregistered}, $\kappa(\widehat{G}(v,1,2))-\kappa(\widetilde{G}(v,1,2))\approx 0.1667$. Further, it is shown in Example \ref{Ex:asymp for star} that $G$ is $(v,1,2)$-paradoxical.}
\end{example}

Let us continue with the hypothesis and notation of Definition \ref{Def:paradoxical}. To see if $G$ is $(v,k_1,k_2)$-paradoxical, we need to investigate $\kappa(\widehat G(v,k_1,k_2))-\kappa(\widetilde G(v,k_1,k_2))$. So, we shall find formulae for  $\mathbf{d}_{\widetilde{G}}^TF_{\widetilde{G}}\mathbf{d}_{\widetilde{G}}$ and $\mathbf{d}_{\widehat G}^TF_{\widehat G}\mathbf{d}_{\widehat G}$ in Lemmas \ref{lem:P_k} and \ref{lem:C_k}, respectively.


\begin{lemma}\label{lem:P_k}
	Let $P_k$ be a path with two pendent vertices $x$ and $y$ where $k\geq 2$, and $H$ be a connected graph. Suppose that $G$ is the graph obtained from $P_k$ and $H$ by identifying a vertex of $P_k$ and a vertex of $H$, say $v$. Let $\mathrm{dist}_G(v,x)=k_1$ and $\mathrm{dist}_G(v,y)=k_2$. Then,
	\begin{align*}
	\hspace{-0.25cm}\mathbf{d}_{G}^TF_{G}\mathbf{d}_{G}=\mathbf{d}_{H}^TF_{H}\mathbf{d}_{H}+4(k-1)\mathbf{d}_H^T\mathbf{f}_H^v+\tau_{H}\left(\frac{4}{3}(k-1)^3+\frac{2}{3}(k-1)+4m_{H}(k_1^2+k_2^2)\right).
	\end{align*}
\end{lemma}
\begin{proof}
	The conclusion is straightforward from \eqref{dfd:expression} and Example \ref{ex:P_n}.
\end{proof}

\begin{lemma}\label{lem:C_k}
	Let $C_k$ be a cycle of length $k$ where $k\geq 3$, and $H$ be a connected graph. Suppose that $G$ is the graph obtained from $C_k$ and $H$ by identifying a vertex of $C_k$ and a vertex of $H$, say $v$. Then,
	\begin{align*}
	\mathbf{d}_{G}^TF_{G}\mathbf{d}_{G}=k\mathbf{d}_H^TF_H\mathbf{d}_H+4k^2\mathbf{d}_{H}^T\mathbf{f}_{H}^v+\frac{2\tau_H}{3}(k+2m_H)(k-1)k(k+1).
	\end{align*}
\end{lemma}
\begin{proof}
	The conclusion is readily established from \eqref{dfd:expression} and Example \ref{ex:C_n}.
\end{proof}


Using Lemmas \ref{lem:P_k} and \ref{lem:C_k}, we establish our desired equivalent condition as follows.

\begin{theorem}\label{Theorem:brass for general twin pendent path}
	Let $G$ be a connected graph with a vertex $v$. Suppose that $k_1,k_2\geq0$, $k_1+k_2\geq 2$ and $k-1=k_1+k_2$. Then, $G$ is $(v,k_1,k_2)$-paradoxical if and only if
	\begin{align}\label{inequality: paradoxical if and only if >0}
	\hspace{-0.3cm}\begin{split}
	&k\mathbf{d}_G^T(2\mathbf{f}_G^v\mathbf{1}^T-F_G)\mathbf{d}_G+4m_G^2\tau_Gk\left(-\frac{2}{3}(k_1+k_2)(k_1+k_2-1)+2k_1k_2)\right)\\
	+&\frac{2m_G\tau_Gk}{3}\left(-5(k_1+k_2)^3+(k_1+k_2)^2+(k_1+k_2)+12k_1k_2(k_1+k_2+1)\right)\\
	-&\frac{2\tau_Gk}{3}(k_1+k_2+1)(k_1+k_2)(k_1+k_2-1)^2>0.
	\end{split}
	\end{align}
\end{theorem}
\begin{proof}
	Evidently, $m_{\widetilde{G}}=m_G+k-1$, $m_{\widehat{G}}=m_G+k$ and $\tau_{\widetilde{G}}=\tau_G$. Since $v$ is a cut-vertex in $\widehat{G}$, we have $\tau_{\widehat{G}}=k\tau_G$. Then,
	\begin{align}\nonumber
	\kappa(\widehat G(v,k_1,k_2))-\kappa(\widetilde G(v,k_1,k_2))&=\frac{\mathbf{d}_{\widehat G}^TF_{\widehat G}\mathbf{d}_{\widehat G}}{4m_{\widehat{G}}\tau_{\widehat{G}}}-\frac{\mathbf{d}_{\widetilde{G}}^TF_{\widetilde{G}}\mathbf{d}_{\widetilde{G}}}{4m_{\widetilde{G}}\tau_{\widetilde{G}}}\\\label{diff:Ghat Gtil}
	&=\frac{(m_G+k-1)\mathbf{d}_{\widehat G}^TF_{\widehat G}\mathbf{d}_{\widehat G}-k(m_G+k)\mathbf{d}_{\widetilde{G}}^TF_{\widetilde{G}}\mathbf{d}_{\widetilde{G}}}{4k(m_G+k)(m_G+k-1)\tau_G}.
	\end{align}
	Then, $G$ is $(v,k_1,k_2)$-paradoxical if and only if $$(m_G+k-1)\mathbf{d}_{\widehat G}^TF_{\widehat G}\mathbf{d}_{\widehat G}-k(m_G+k)\mathbf{d}_{\widetilde{G}}^TF_{\widetilde{G}}\mathbf{d}_{\widetilde{G}}>0.$$ 
	For simplicity, let $\mathbf{d}=\mathbf{d}_G$, $\mathbf{f}^v=\mathbf{f}_G^v$, $F=F_G$, $m=m_G$ and $\tau=\tau_G$. Using Lemmas \ref{lem:P_k} and \ref{lem:C_k}, we have 
	\begin{align}\label{temp:eqn11}
	\begin{split}
	&(m+k-1)\mathbf{d}_{\widehat G}^TF_{\widehat G}\mathbf{d}_{\widehat G}-k(m+k)\mathbf{d}_{\widetilde{G}}^TF_{\widetilde{G}}\mathbf{d}_{\widetilde{G}}\\
	=&(m+k-1)\left(k\mathbf{d}^TF\mathbf{d}+4k^2\mathbf{d}^T\mathbf{f}^v+\frac{2\tau}{3}(k+2m)(k-1)k(k+1)\right)\\
	&-k(m+k)\left(\mathbf{d}^TF\mathbf{d}+4(k-1)\mathbf{d}^T\mathbf{f}^v+\frac{4}{3}\tau(k-1)^3+\frac{2}{3}\tau(k-1)+4m\tau(k_1^2+k_2^2)\right)\\
	=&-k\mathbf{d}^TF\mathbf{d}+4mk\mathbf{d}^T\mathbf{f}^v+4m^2\tau k\left(\frac{1}{3}(k-1)(k+1)-k_1^2-k_2^2\right)\\
	+&\frac{2m\tau k}{3}\left((k-1)k(k+1)+2(k-1)^2(k+1)-2(k-1)^3-(k-1)-6k(k_1^2+k_2^2)\right)\\
	+&\frac{2\tau k}{3}\left((k-1)^2k(k+1)-2k(k-1)^3-k(k-1)\right).
	\end{split}
	\end{align}
	Since $\mathbf{1}^T\mathbf{d}=2m$, we have $4mk\mathbf{d}^T\mathbf{f}^v=2k\mathbf{d}^T\mathbf{f}^v\mathbf{1}^T\mathbf{d}$. Then, one can check from $k-1=k_1+k_2$ that the last expression in \eqref{temp:eqn11} can be recast as the left side of the inequality \eqref{inequality: paradoxical if and only if >0}.
\end{proof}

Now, we shall introduce the following notation to easily analyse the expression in \eqref{inequality: paradoxical if and only if >0}. Let $G$ be a connected graph of order $n$ with $V(G)=\{1,\dots,n\}$, and let
\begin{align*}
\phi_G(v)&:=\mathbf{d}_G^T(2\mathbf{f}_G^v\mathbf{1}^T-F_G)\mathbf{d}_G,\\
\phi_1(k_1,k_2)&:=-\frac{2}{3}(k_1+k_2)(k_1+k_2-1)+2k_1k_2,\\
\phi_2(k_1,k_2)&:=-(k_1+k_2)(5(k_1+k_2)^2-(k_1+k_2)-1)+12k_1k_2(k_1+k_2+1),\\
\phi_3(k_1,k_2)&:=-(k_1+k_2+1)(k_1+k_2)(k_1+k_2-1)^2,
\end{align*}
where $v$, $k_1$ and $k_2$ are integers such that $1\leq v\leq n$, $k_1, k_2\geq 0$ and $k_1+k_2\geq 2$. Furthermore, let
\begin{align}\label{Formula: Phi_G}
\hspace{-0.1cm}\Phi_G(v,k_1,k_2):=k\phi_G(v)+4m_G^2\tau_Gk\phi_1(k_1,k_2)+\frac{2m_G\tau_Gk}{3}\phi_2(k_1,k_2)+\frac{2\tau_G k}{3}\phi_3(k_1,k_2).
\end{align} 
By Theorem \ref{Theorem:brass for general twin pendent path}, $G$ is $(v,k_1,k_2)$-paradoxical if and only if $\Phi_G(v,k_1,k_2)>0$. We simply write $\Phi_G(v,k_1,k_2)$ and $\phi_G(v)$ as $\Phi(v,k_1,k_2)$ and $\phi(v)$, respectively, if $G$ is clear from the context. Note that $\phi_i(k_1,k_2)=\phi_i(k_2,k_1)$ for $i=1,2,3$. So, $\Phi_G(v,k_1,k_2)=\Phi_G(v,k_2,k_1)$.
\begin{remark}
	{\rm A connected graph $G$ is $(v,k_1,k_2)$-paradoxical if and only if $G$ is $(v,k_2,k_1)$-paradoxical. Furthermore, $G$ is $(k_1,k_2)$-paradoxical if and only if $G$ is $(k_2,k_1)$-paradoxical.}
\end{remark}

\subsection*{Signs of $\phi_i(k_1,k_2)$ for $i=1,2,3$ and $\phi_G(v)$}
We shall consider the signs of $\phi_i(k_1,k_2)$ for $i=1,2,3$ in terms of $k_1$ and $k_2$, and consider an upper bound for each $\phi_i(k_1,k_2)$. Evidently, $\phi_3(k_1,k_2)$ decreases as $k_1+k_2$ increases, and so
\begin{align}\label{minimum of phi_3}
\phi_3(k_1,k_2)\leq -6\;\text{for any $k_1, k_2\geq 0$ where $k_1+k_2\geq 2$}
\end{align}
with equality if and only if $k_1+k_2=2$. Next, $\phi_1(k_1,k_2)$ can be written as 
$$\phi_1(k_1,k_2)=-\frac{2}{3}(k_1^2-(k_2+1)k_1+k_2^2-k_2).$$
Setting $\phi_1(k_1,k_2)=0$, we have 
$$
k_1=\frac{1}{2}\left((k_2+1)\pm\sqrt{-3k_2^2+6k_2+1}\right).
$$
Since $\phi_1(k_1,k_2)$ is symmetric, without loss of generality, we shall fix $k_2$ first. It follows from $-3k_2^2+6k_2+1<0$ that if $k_2<1-\frac{2\sqrt{3}}{3}<0$ or $k_2>1+\frac{2\sqrt{3}}{3}>2$, then $\phi_1(k_1,k_2)<0$ for any $k_1\geq 0$. Furthermore, if $k_2=1$, then $\phi_1(1,1)=\frac{2}{3}$, $\phi_1(2,1)=0$ and $\phi_1(k_1,1)<0$ for $k_1>2$. Finally, for $k_2=2$, we have $\phi_1(0,2)=-\frac{4}{3}$, $\phi_1(1,2)=\phi_1(2,2)=0$ and $\phi_1(k_1,2)<0$ for $k_1>2$. Therefore, $\phi_1(k_1,k_2)>0$ if and only if $(k_1,k_2)=(1,1)$; $\phi_1(k_1,k_2)=0$ if and only if $(k_1,k_2)\in\{(1,2),(2,1),(2,2)\}$; and $\phi_1(k_1,k_2)<0$ for any $k_1,k_2\geq 0$ with $k_1+k_2\geq 2$ and $(k_1,k_2)\notin\{(1,1),(1,2),(2,1),(2,2)\}$.

\begin{remark}\label{Remark:phi_1 decrease}
	{\rm We have $\frac{\partial \phi_1}{\partial k_1}=-\frac{4}{3}k_1+\frac{2}{3}(k_2+1)$. Then, $\phi_1(2,0)=-\frac{4}{3}$ and $\frac{\partial \phi_1}{\partial k_1}\Bigr|_{\substack{k_2=0}}<0$ for $k_1\geq 2$; $\phi_1(3,1)=-2$ and $\frac{\partial \phi_1}{\partial k_1}\Bigr|_{\substack{k_2=1}}<0$ for $k_1\geq 3$; $\phi_1(3,2)=-\frac{4}{3}$ and $\frac{\partial \phi_1}{\partial k_1}\Bigr|_{\substack{k_2=2}}<0$ for $k_1\geq 3$; finally, $\phi_1(k_2,k_2)=-\frac{2}{3}(k_2^2-2k_2)\leq -2$ for $k_2\geq 3$ and  $\frac{\partial \phi_1}{\partial k_1}<0$ for $k_1\geq k_2\geq 3$. Hence, since $\phi_1(k_1,k_2)$ is symmetric, $\phi_1(k_1,k_2)\leq -\frac{4}{3}$ for integers $k_1,k_2\geq 0$ with $k_1+k_2\geq 2$ and $(k_1,k_2)\notin\{(1,1),(1,2),(2,1),(2,2)\}$. Furthermore, by computation, we have $\phi_1(3,0)=\phi_1(4,2)=-4$. Therefore, $\phi_1(k_1,k_2)\leq -2$ for integers $k_1,k_2\geq 0$ with $k_1+k_2\geq 2$ and $(k_1,k_2)\notin\{(0,2),(2,0),(1,1),(1,2),(2,1),(2,2),(2,3),(3,2)\}$.}
\end{remark}

Putting $k-1=k_1+k_2\geq 2$, $\phi_2(k_1,k_2)$ can be written as
$$
\phi_2(k_1,k_2)=-12kk_1^2+12k(k-1)k_1-(k-1)(5k^2-11k+5).
$$
Setting $\phi_2(k_1,k_2)=0$, we have
$$
k_1=\frac{1}{12k}\left(6k(k-1)\pm\sqrt{-12k(k-1)(2k^2-8k+5)}\right).
$$
Since $2k^2-8k+5>0$ for all $k\geq 4$, $\phi_2(k_1,k_2)<0$ for any $k_1,k_2\geq 0$ with $k_1+k_2\geq 3$. For $k=3$, we have $\phi_2(1,1)=2>0$ and $\phi_2(2,0)=-34<0$. Let $f(t)=-(t-1)(2t^2-8t+5)$ where $t$ is real number. Then, for fixed $t\geq 3$, the maximum of $\phi_2(t_1,t_2)$ for nonnegative numbers $t_1$ and $t_2$ with $t_1+t_2=t-1$ is attained as $f(t)$ at $t_1=\frac{t-1}{2}$. We can find that $f(3)>0$, $f(4)=-15$ and $f'(t)<0$ for $t\geq 4$. From computation, we have $\phi_2(0,3)=-123$ and $\phi_2(1,2)=-27$. Hence,
\begin{align}\label{minimum of phi_2 for k}
\phi_2(k_1,k_2)<-15\; \text{for any $k_1,k_2\geq 0$ with $k_1+k_2\geq 2$ and $(k_1,k_2)\neq(1,1)$.}
\end{align}

We claim that for a non-trivial connected graph $G$, $\phi(v)=\mathbf{d}^T(2\mathbf{f}^v\mathbf{1}^T-F)\mathbf{d}>0$ for $v=1,\dots,n$. In order to establish our claim, we first show that $f_{i,j}^G$ is a metric on the vertex set of $G$ by using the resistance distance (see \cite{Klein:resistance} for an introduction). Let $L$ be the Laplacian matrix of $G$, and let $L^\dagger=[\ell^\dagger_{i,j}]_{n\times n}$ be the Moore--Penrose inverse of $L$. Then, the resistance distance $\Omega_{i,j}$ between vertices $i$ and $j$ of $G$ is represented (see \cite{Klein:sum}) as:
$$
\Omega_{i,j}=\ell^\dagger_{i,i}+\ell^\dagger_{j,j}-\ell^\dagger_{i,j}-\ell^\dagger_{j,i}.
$$
Moreover, the resistance distance is a metric on $V(G)$ (see \cite{Bapat:graphs}). As proved in \cite{Pavel:resistance}, the number $f_{i,j}^G$ of $2$-tree spanning forests of $G$ having $i$ and $j$ in different trees is
$$
f_{i,j}^G=\tau_G\Omega_{i,j}.
$$
Therefore, we have the following properties endowed by the metric $\Omega_{i,j}$:
\begin{enumerate}[label=(\roman*)]
	\item $f_{i,j}^G\geq 0$, with equality if and only if $i=j$;
	\item $f_{i,j}^G=f_{j,i}^G$ for all $i$, $j$;
	\item for any $i,j,k$, $f_{i,j}^G\leq f_{i,k}^G+f_{k,j}^G$, with equality \cite{Kirkland:CombinatorialBook} if and only if either all paths in $G$ from $i$ to $j$ pass through $k$ or $k$ is one of $i$ and $j$.
\end{enumerate}
Let $X=2\mathbf{f}^v\mathbf{1}^T-F$, and $Q=[q_{i,j}]=\frac{X+X^T}{4}$. Then,
$$\mathbf{d}^TX\mathbf{d}=\frac{1}{2}(\mathbf{d}^TX\mathbf{d}+\mathbf{d}^TX^T\mathbf{d})=2\mathbf{d}^TQ\mathbf{d}.$$ 
Since $2Q=\mathbf{f}^v\mathbf{1}^T+\mathbf{1}(\mathbf{f}^v)^T-F$, we have $2q_{i,j}=f_{i,v}+f_{v,j}-f_{i,j}\geq 0$. Since $G$ is connected, if $i\neq v$, then there exists a $2$-tree spanning forest having $i$ and $v$ in different trees, \textit{i.e.,} $f_{i,v}>0$. For a non-trivial connected graph $G$, there exists a vertex $i$ with $i\neq v$ such that $2q_{i,i}=f_{i,v}+f_{v,i}>0$. Hence, $Q$ is a non-negative symmetric matrix with $Q\neq 0$. Since $\mathbf{d}>0$, we have $\mathbf{d}^TQ\mathbf{d}>0$. Therefore, $\mathbf{d}^T(2\mathbf{f}^v\mathbf{1}^T-F)\mathbf{d}>0$.

Recall that given a connected graph $G$ of order $n$ where $n\geq 2$, $G$ is $(v,k_1,k_2)$-paradoxical if and only if $\Phi(v,k_1,k_2)>0$, where $1\leq v\leq n$ and integers $k_1,k_2\geq 0$ with $k_1+k_2\geq 2$. We have seen that $\phi(v)>0$ for any $1\leq v\leq n$ regardless of $k_1$ and $k_2$; $\phi_1(k_1,k_2)\geq 0$ if and only if $(k_1,k_2)\in\{(1,1),(1,2),(2,1),(2,2)\}$; $\phi_2(k_1,k_2)<0$ for any $k_1,k_2$ with $(k_1,k_2)\neq(1,1)$; and $\phi_3(k_1,k_2)<0$ for any $k_1,k_2$. Hence, $\phi(v)$ must have a relatively larger quantity in order for $G$ to be $(v,k_1,k_2)$-paradoxical.

Consider the case $k_1=k_2=1$. Then, $\Phi(v,1,1)=3\phi(v)+8m_G^2\tau_G+4m_G\tau_G-12\tau_G$. Clearly, $\Phi(v,1,1)>0$ for any non-trivial connected graph $G$ and any vertex $v$ of $G$. Hence, we have the following result.

\begin{theorem}\cite{CiardoparadoxicalTwins}\label{Theorem: (v,1,1)-paradoxical}
	Let $G$ be a connected graph of order $n$ where $n\geq 2$. Then, $G$ is $(1,1)$-paradoxical.
\end{theorem}

\subsection*{Combinatorial interpretation for $\mathbf{f}^v_G\mathbf{1}^T+\mathbf{1}(\mathbf{f}^v_G)^T-F_G$}
Let $G$ be a non-trivial connected graph with a vertex $v$. We now discuss a combinatorial interpretation for $q_{i,j}$ where $q_{i,j}=\frac{1}{2}(f_{i,v}^G+f_{v,j}^G-f_{i,j}^G)$. Denote by $\mathcal{F}_{G}(i,j;v)$ (or equivalently $\mathcal{F}_{G}(v;i,j)$) the set of all spanning forests consisting of two trees in $G$, one of which contains vertices $i$ and $j$ and the other of which contains a vertex $v$. Then, we have
\begin{align*}
|\mathcal{F}_{G}(i;j)|&=|\mathcal{F}_{G}(i;v,j)|+|\mathcal{F}_{G}(i,v;j)|,\\
|\mathcal{F}_{G}(i;v)|&=|\mathcal{F}_{G}(i,j;v)|+|\mathcal{F}_{G}(i;v,j)|,\\
|\mathcal{F}_{G}(v;j)|&=|\mathcal{F}_{G}(i,v;j)|+|\mathcal{F}_{G}(v;i,j)|.
\end{align*}
It follows that $2q_{i,j}=f_{i,v}+f_{v,j}-f_{i,j}=2|\mathcal{F}_{G}(i,j;v)|$, that is, $q_{i,j}$ is the number of $2$-tree spanning forests of $G$ having $i$, $j$ in one tree and $v$ in the other. Thus, we define $Q_{G,v}$ as the matrix $Q_{G,v}=[q_{i,j}]$ associated to $G$ and $v$. Then,
\begin{align*}
Q_{G,v}=\frac{1}{2}(\mathbf{f}^v\mathbf{1}^T+\mathbf{1}(\mathbf{f}^v)^T-F),\;\;\phi_G(v)=\mathbf{d}^T(2\mathbf{f}^v\mathbf{1}^T-F)\mathbf{d}=2\mathbf{d}^TQ_{G,v}\mathbf{d}.
\end{align*}

\begin{remark}\label{Remark:entries r}
	{\rm Let $G$ be a connected graph with a vertex $v$. Let $Q_{G,v}=[q_{i,j}]$. Since $2q_{i,j}=f_{i,v}+f_{v,j}-f_{i,j}$, we have $q_{i,j}=0$ whenever $v=i$ or $v=j$. Suppose that $v$ is a cut-vertex. If there is no path from $i$ to $j$ with $i\neq v$ and $j\neq v$ in $G-v$, then by the combinatorial interpretation for $q_{i,j}$, we obtain $q_{i,j}=0$. Consider a branch $B$ of $G$ at $v$. Let $i,j\in V(B)$. For each forest in $\mathcal{F}_{G}(i,j;v)$, the subtree with the vertex $v$ in the forest must contain all vertices of $V(G)\backslash V(B)$. Thus, for the subgraph $G'$ induced by $V(G)\backslash (V(B)-\{v\})$, we have $|\mathcal{F}_{G}(i,j;v)|=\tau_{G'}|\mathcal{F}_{B}(i,j;v)|$. This implies that $G'$ is a tree if and only if $|\mathcal{F}_{G}(i,j;v)|=|\mathcal{F}_{B}(i,j;v)|$.
		
		Given a tree $\cT$ with a vertex $v$, let $Q_{\cT,v}=[q_{i,j}]$. Consider two vertices $i$ and $j$ in $\cT$ with $i\neq v$ and $j\neq v$. For each forest in $\mathcal{F}_{\cT}(i,j;v)$, there is a subtree of the forest having $i$ and $j$. Then, all vertices $w_0,w_1,\dots,w_{\mathrm{dist}(i,j)}$ on the subpath with pendent vertices $i$ and $j$ must be contained in the subtree. Therefore, $q_{i,j}=\mathrm{min}\{\mathrm{dist}(v,w_p)|p=0,\dots,\mathrm{dist}(i,j)\}$. In particular, if $i=j$ then $q_{i,j}=\mathrm{dist}(i,v)$.}
\end{remark}

Based on Remark \ref{Remark:entries r}, let us consider the following example.

\begin{example}\label{Example:R_Pn v}
	{\rm Consider the path $P_6=(1,\dots,6)$. Let $Q_{P_n,v}=[q_{i,j}]$ where $v=3$. Evidently, $q_{3,i}=q_{i,3}=0$ for $1\leq i\leq 6$. Since $v$ is a cut vertex, we have $q_{i,j}=0$ for $i\in\{1,2\}$ and $j\in\{4,5,6\}$. By the argument in the second paragraph of Remark \ref{Remark:entries r}, we have
		\begin{align*}
		Q_{P_n,v}=\begin{bmatrix}
		2 & 1 & 0 & 0 & 0 & 0\\
		1 & 1 & 0 & 0 & 0 & 0\\
		0 & 0 & 0 & 0 & 0 & 0\\
		0 & 0 & 0 & 1 & 1 & 1\\
		0 & 0 & 0 & 1 & 2 & 2\\
		0 & 0 & 0 & 1 & 3 & 3\\
		\end{bmatrix}.
		\end{align*}}
\end{example}

While the combinatorial interpretation for entries of $Q_{G,v}$ is given, we mainly focus on the computation of $\mathbf{d}^TQ_{G,v}\mathbf{d}$ in the rest of this section; but the combinatorial interpretation is used more in Section \ref{Section:Asymptotic for trees}.

\subsection*{Several examples} We now find conditions for $K_n$, $C_n$, $P_n$ or $S_n$ to be $(v,k_1,k_2)$-paradoxical or $(k_1,k_2)$-paradoxical. For simplicity, set $k-1=k_1+k_2$ and $\phi_i=\phi_i(k_1,k_2)$ for $i=1,2,3$. Note that $\phi_G(v)=\mathbf{d}^T(2\mathbf{f}^v\mathbf{1}^T-F)\mathbf{d}=2\mathbf{d}^TQ_{G,v}\mathbf{d}$ and $\phi_i(k_1,k_2)=\phi_i(k_2,k_1)$ for $i=1,2,3$. We compute $\phi_G(v)$ by using $F$ and $\mathbf{f}^v$ for $G=K_n$ or $G=C_n$, and by directly finding $Q_{G,v}$ for $G=P_n$ or $G=S_n$. For convenience, the following quantities are computed in advance: $\phi_2(1,2)=-27$, $\phi_3(1,2)=-48$, $\phi_2(2,2)=-60$ and $\phi_3(2,2)=-180$.

\begin{example}\label{Ex:asymp for K_n}
	{\rm Consider a complete graph $K_n$ on $n$ vertices. Let $v$ be a vertex of $K_n$. Then, from $\phi(v)=\mathbf{d}^T(2\mathbf{f}^v\mathbf{1}^T-F)\mathbf{d}$ and Example \ref{ex:K_n}, it is readily seen that
		$$
		\phi(v)=2n^{n-2}(n-1)^3=2\tau(n-1)^3.
		$$ 
		Using \eqref{Formula: Phi_G} with $m=\frac{n(n-1)}{2}$, we obtain
		\begin{align*}
		\Phi_{K_n}(v,k_1,k_2)=\tau k\left(2(n-1)^3+n^2(n-1)^2\phi_1+\frac{1}{3}n(n-1)\phi_2+\frac{2}{3}\phi_3\right).
		\end{align*}
		Suppose that $(k_1,k_2)\notin \{(1,1),(1,2),(2,1),(2,2)\}$. By Remark \ref{Remark:phi_1 decrease}, $\phi_1\leq -\frac{4}{3}$. From \eqref{minimum of phi_2 for k}, we have $\phi_2<-15$. By \eqref{minimum of phi_3}, $\phi_3\leq -6$. Hence,
		\begin{align*}
		\Phi(v,k_1,k_2)<\tau k\left(-\frac{4}{3}n^4 + \frac{14}{3}n^3-\frac{37}{3}n^2 +11n-6\right).
		\end{align*}
		One can verify that $-\frac{4}{3}n^4 + \frac{14}{3}n^3-\frac{37}{3}n^2 +11n-6<0$ for $n\geq 1$. Thus, if $(k_1,k_2)\notin\{(1,1),(1,2),(2,1),(2,2)\}$, then $K_n$ is not $(v,k_1,k_2)$-paradoxical for any $n\geq 1$. 
		
		Consider $(k_1,k_2)=(1,2)$ and $(k_1,k_2)=(2,2)$. Then,
		\begin{align*}
		&\Phi(v,1,2)=4\tau\left(2(n-1)^3-9n(n-1)-32\right),\\
		&\Phi(v,2,2)=5\tau\left(2(n-1)^3-20n(n-1)-120\right).
		\end{align*}
		Using the derivatives of $\frac{\Phi(v,1,2)}{4\tau}$ and $\frac{\Phi(v,2,2)}{10\tau}$ with respect to $n$, it can be checked that $\Phi(v,1,2)>0$ if and only if $n\geq 7$; $\Phi(v,2,2)>0$ if and only if $n\geq 13$. Hence, $K_n$ is $(1,2)$-paradoxical for $n\geq 7$, and $(2,2)$-paradoxical for $n\geq 13$.}
\end{example}

\begin{example}\label{Ex:asymp for C_n}
	{\rm Given a cycle $C_n$ with a vertex $v$, from $\phi(v)=\mathbf{d}^T(2\mathbf{f}^v\mathbf{1}^T-F)\mathbf{d}$ and Example \ref{ex:C_n}, we have $\phi(v)=\frac{2}{3}(n-1)n^2(n+1)=\frac{2\tau}{3}(n-1)n(n+1)$. Using \eqref{Formula: Phi_G}, we find
		\begin{align*}
		\Phi_{C_n}(v,k_1,k_2)=\tau k\left(\frac{2}{3}(n-1)n(n+1)+4n^2\phi_1+\frac{2}{3}n\phi_2+\frac{2}{3}\phi_3\right).
		\end{align*}
		We observe that the term of the highest degree about $n$ in $\frac{\Phi(v,k_1,k_2)}{\tau k}$ has a positive coefficient. This implies that given $k_1,k_2\geq 0$ with $k_1+k_2\geq 2$, $C_n$ is $(k_1,k_2)$-paradoxical for sufficiently large $n$. Consider $(k_1,k_2)=(1,2)$ and $(k_1,k_2)=(2,2)$. Then,
		\begin{align*}
		&\Phi(v,1,2)=4\tau\left(\frac{2}{3}(n-1)n(n+1)-18n-32\right),\\
		&\Phi(v,2,2)=5\tau\left(\frac{2}{3}(n-1)n(n+1)-40n-120\right).
		\end{align*}
		One can verify that $\Phi(v,1,2)\geq 0$ for $n\geq 6$ with equality if and only if  $n=6$; $\Phi(v,2,2)\geq 0$ for $n\geq 9$ with equality if and only if $n=9$. Hence, $C_n$ is $(1,2)$-paradoxical for $n\geq 7$, and $(2,2)$-paradoxical for $n\geq 10$.}
\end{example}

\begin{example}\label{Ex:asymp for P_n}
	{\rm Consider the path $P_n=(1,\dots,n)$ with a vertex $v$. By Remark \ref{Remark:entries r} and Example \ref{Example:R_Pn v}, we have
		\begin{align*}
		Q_{P_n,v}=\begin{bmatrix}
		M_1 & 0 \\
		0 & M_2
		\end{bmatrix}
		\end{align*}
		where $M_1=\begin{bmatrix}
		\mathrm{min}\{v-i,v-j\}
		\end{bmatrix}_{1\leq i,j\leq v}$ and $M_2=\begin{bmatrix}
		\mathrm{min}\{i,j\}
		\end{bmatrix}_{1\leq i,j\leq n-v}$. We have $\mathbf{d}_{P_n}=2\mathbf{1}_{n}-\mathbf{e}_1-\mathbf{e}_{n}$. Then,
		\begin{align*}
		\mathbf{d}^TQ_{P_n,v}\mathbf{d}&=4\mathbf{1}^TQ_{P_n,v}\mathbf{1}+(M_1)_{1,1}+(M_2)_{n-v,n-v}-4\mathbf{1}^TQ_{P_n,v}\be_1-4\mathbf{1}^TQ_{P_n,v}\be_n\\
		&=4\left(\sum_{k=1}^{v-1}k^2+\sum_{k=1}^{n-v}k^2\right)+n-1-2v(v-1)-2(n-v)(n-v+1)\\
		&=4(n-1)v^2-4(n^2-1)v+\frac{4}{3}n^3-\frac{1}{3}n-1.
		\end{align*}
		The minimum of $\mathbf{d}^TQ_{P_n,v}\mathbf{d}$ is attained as $\frac{1}{3}n(n-1)(n-2)$ if $n$ is odd, and as $\frac{1}{3}n^3-n^2+\frac{5}{3}n-1$ if $n$ is even. The maximum of $\mathbf{d}^TQ_{P_n,v}\mathbf{d}$ is $\frac{1}{3}(n-1)(2n-1)(2n-3)$ at $v=1$ or $v=n$. 
		
		By \eqref{Formula: Phi_G} and the minimum of $\phi(v)=2\mathbf{d}^TQ_{P_n,v}\mathbf{d}$, we have
		$$\Phi_{P_n}(v,k_1,k_2)\geq k\left(\frac{2}{3}n(n-1)(n-2)+4(n-1)^2\phi_1+\frac{2}{3}(n-1)\phi_2+\frac{2}{3}\phi_3\right).$$
		By a similar argument as in Example \ref{Ex:asymp for C_n}, given $k_1,k_2\geq 0$ with $k_1+k_2\geq 2$, $P_n$ is $(k_1,k_2)$-paradoxical for sufficiently large $n$.}  
\end{example}

\begin{example}\label{Ex:asymp for star}
	{\rm Consider a star $S_n$ of order $n$ with a vertex $v$. Using Remark \ref{Remark:entries r}, it can be checked that
		\begin{align*}
		Q_{S_n,v}=\begin{cases*}
		J+\begin{bmatrix}
		I_{n-1} & 0\\
		0 & 0
		\end{bmatrix}-\be_v\mathbf{1}^T-\mathbf{1}\be_v^T, & \text{if $\mathrm{deg}(v)=1$,}\\
		\begin{bmatrix}
		I_{n-1} & 0\\
		0 & 0
		\end{bmatrix}, & \text{if $\mathrm{deg}(v)=n-1$}.
		\end{cases*}
		\end{align*}
		Hence,
		$$\mathbf{d}^TQ_{S_n,v}\mathbf{d}=\begin{cases*}
		(n-1)(4n-7), & \text{if $\mathrm{deg}(v)=1$,}\\
		n-1, & \text{if $\mathrm{deg}(v)=n-1$}.
		\end{cases*}$$ 
		
		Suppose that $v$ is the centre vertex. Then, $n\geq 3$. By \eqref{Formula: Phi_G} and $\phi(v)=2\mathbf{d}^TQ_{S_n,v}\mathbf{d}$, we have
		\begin{align*}
		\Phi_{S_n}(v,k_1,k_2)=k\left(2(n-1)+4(n-1)^2\phi_1+\frac{2}{3}(n-1)\phi_2+\frac{2}{3}\phi_3\right).
		\end{align*}
		Let $(k_1,k_2)\neq (1,1)$. Clearly, $\phi_1\leq 0$. By \eqref{minimum of phi_2 for k} and \eqref{minimum of phi_3}, we have $\phi_2<-15$ and $\phi_3\leq -6$, respectively. So, $\Phi(v,k_1,k_2)<-4k(2n-1)$. Hence, if $S_n$ is $(v,k_1,k_2)$-paradoxical where $v$ is the centre vertex of $S_n$, then $(k_1,k_2)=(1,1)$ and $n\geq 3$. 
		
		Suppose that  $v$ is a pendent vertex. Then,
		\begin{align*}
		\Phi_{S_n}(v,k_1,k_2)=k\left(2(n-1)(4n-7)+4(n-1)^2\phi_1+\frac{2}{3}(n-1)\phi_2+\frac{2}{3}\phi_3\right).
		\end{align*}
		We have $\phi_1(2,0)=-\frac{4}{3}$, $\phi_2(2,0)=-34$ and $\phi_3(2,0)=-6$; $\phi_1(3,2)=-\frac{4}{3}$, $\phi_2(3,2)=-163$ and $\phi_3(3,2)=-480$. One can check that $\Phi(v,2,0)=8n^2-102n+82>0$ for $n\geq 12$; $\Phi(v,2,1)=32n^2-160n>0$ for $n\geq 6$; $\Phi(v,2,2)=40n^2-310n-330>0$ for $n\geq 9$; and $\Phi(v,3,2)=16n^2-720n-1216>0$ for $n\geq 47$. Let $$A=\{(0,2),(2,0),(1,1),(1,2),(2,1),(2,2),(2,3),(3,2)\}.$$
		Suppose that $(k_1,k_2)\notin A$. By Remark \ref{Remark:phi_1 decrease}, we have $\phi_1(k_1,k_2)\leq -2$. From \eqref{minimum of phi_2 for k} and \eqref{minimum of phi_3}, $\phi_2< -15$ and $\phi_3\leq -6$, respectively. Hence, $\Phi(v,k_1,k_2)< -k(16n-12)$. Therefore, if $S_n$ is $(v,k_1,k_2)$-paradoxical for a pendent vertex $v$, then $k_1$, $k_2$ and $n$ satisfy one of the following: \begin{enumerate*}[label=(\roman*)]
			\item $(k_1,k_2)=(1,1)$, $n\geq 2$; \item $(k_1,k_2)\in\{(0,2),(2,0)\}$, $n\geq 12$; \item $(k_1,k_2)\in\{(1,2),(2,1)\}$, $n\geq 6$; \item $(k_1,k_2)=(2,2)$, $n\geq 9$; and \item $(k_1,k_2)\in\{(2,3),(3,2)\}$, $n\geq 47$.
	\end{enumerate*}} 
\end{example}

\section{Asymptotic behaviour of a sequence of graphs with twin pendent paths regarding the Braess edge}\label{Section: Asymptotic}

We have seen the families of complete graphs, cycles, stars, and paths in the previous section, and we have observed their asymptotic behaviours with respect to the property of being $(v,k_1,k_2)$-paradoxical as the orders of graphs increase. In particular, from Examples \ref{Ex:asymp for C_n} and \ref{Ex:asymp for P_n}, if for any non-negative integers $k_1$ and $k_2$ with $k_1+k_2\geq 2$, any graph in a family of cycles or paths has sufficiently large order relative to $k_1$ and $k_2$, then it is $(k_1,k_2)$-paradoxical. This idea is formalized for a specified vertex, and a tool for finding such families is described in this section.

\begin{definition}\label{Def:asymptotically}
	{\rm Let $\cG^v$ be a sequence of graphs $G_1,G_2,\dots$ where for each $n\geq 1$, $G_n$ is a connected graph of order $n$ with a {\em specified} vertex $v$. Fix integers $k_1,k_2\geq 0$ with $k_1+k_2\geq 2$. The sequence $\cG^v$ is {\em asymptotically $(k_1,k_2)$-paradoxical} if there exists $N>0$ such that $G_n$ is $(v,k_1,k_2)$-paradoxical for all $n\geq N$. The sequence $\cG^v$ is {\em asymptotically paradoxical} if for any integers $l_1,l_2\geq 0$ with $l_1+l_2\geq 2$, $\cG^v$ is asymptotically $(l_1,l_2)$-paradoxical.}
\end{definition}


In what follows, $\cG^v=(G_n)^v$ denotes a sequence of connected graphs $G_1,G_2,\dots$ where for each $n\geq 1$, $|V(G_n)|=n$ and $v\in V(G_n)$.

\begin{example}
	{\rm From Theorem \ref{Theorem: (v,1,1)-paradoxical}, any sequence $\cG^v=(G_n)^v$ is asymptotically $(1,1)$-paradoxical.}
\end{example}

\begin{example}
	{\rm Let $\cG_1^v=(K_n)^v$, $\cG_2^v=(C_n)^v$, $\cG_3^v=(P_n)^v$ and $\cG_4^v=(S_n)^v$. From Examples \ref{Ex:asymp for K_n}--\ref{Ex:asymp for star}, $\cG_2^v$ and $\cG_3^v$ are asymptotically paradoxical, but $\cG_1^v$ and $\cG_4^v$ are not. In particular, $\cG_1^v$ is asymptotically $(k_1,k_2)$-paradoxical if and only if $(k_1,k_2)\in\{(1,1),(1,2),(2,1),(2,2)\}$. Consider $\cG_4^v=(S_n)^v$. Suppose that there exists $N>0$ such that $v$ is a pendent vertex of $S_n$ for all $n\geq N$. Then, $\cG_4^v$ is asymptotically $(k_1,k_2)$-paradoxical if and only if $(k_1,k_2)$ is in the set $A$ described in Example \ref{Ex:asymp for star}. If there exists $N>0$ such that $v$ is the centre vertex of $S_n$ for all $n\geq N$, then $\cG_4^v$ is asymptotically $(k_1,k_2)$-paradoxical if and only if $(k_1,k_2)=(1,1)$.}
\end{example}

Consider a sequence $\cG^v=(G_n)^v$. Examining the proof of Theorem \ref{Theorem:brass for general twin pendent path} with \eqref{diff:Ghat Gtil}, we find from \eqref{Formula: Phi_G} that
\begin{align}\nonumber
&\kappa(\widehat{G}_n(v,k_1,k_2))-\kappa(\widetilde{G}_n(v,k_1,k_2))\\\nonumber
=&\frac{\Phi_{G_n}(v,k_1,k_2)}{4k(m_{G_n}+k)(m_{G_n}+k-1)\tau_{G_n}}\\\label{Tempeqn}
=&\frac{\phi_{G_n}(v)+4m_{G_n}^2\tau_{G_n}\phi_1(k_1,k_2)+\frac{2m_{G_n}\tau_{G_n}}{3}\phi_2(k_1,k_2)+\frac{2\tau_{G_n}}{3}\phi_3(k_1,k_2)}{4(m_{G_n}+k)(m_{G_n}+k-1)\tau_{G_n}}.
\end{align}
Note that since $\phi_{G_n}(v)>0$ for all $n\geq 2$, we have $\frac{\phi_{G_n}(v)}{4m_{G_n}^2\tau_{G_n}}>0$.

We introduce a sufficient condition for $\cG^v=(G_n)^v$ to be asymptotically $(k_1,k_2)$-paradoxical. Moreover, the following result can be used for minimizing $N_0>0$ such that $G_n$ is $(v,k_1,k_2)$-paradoxical for all $n\geq N_0$.


\begin{proposition}\label{Proposition:monotone}
	Let $k_1,k_2$ be non-negative integers with $k_1+k_2\geq 2$ and $(k_1,k_2)\neq(1,1)$. Given a sequence $\cG^v=(G_n)^v$, suppose that $\frac{\phi_{G_{N+1}}(v)}{4m_{G_{N+1}}^2\tau_{G_{N+1}}}\geq\frac{\phi_{G_N}(v)}{4m_{G_N}^2\tau_{G_N}}$ for some $N>0$. If $G_N$ is $(v,k_1,k_2)$-paradoxical, then $G_{N+1}$ is $(v,k_1,k_2)$-paradoxical. This implies that if $G_{N_0}$ is $(v,k_1,k_2)$-paradoxical for some $N_0>0$, and if $\frac{\phi_{G_n}(v)}{4m_{G_n}^2\tau_{G_n}}$ is non-decreasing for $n\geq N_0$, then $G_n$ is $(v,k_1,k_2)$-paradoxical for all $n\geq N_0$---that is, $\cG^v=(G_n)^v$ is asymptotically $(k_1,k_2)$-paradoxical.
\end{proposition}
\begin{proof}
	We only need to show that if $\Phi_{G_{N}}(v,k_1,k_2)>0$ then $\Phi_{G_{N+1}}(v,k_1,k_2)>0$. We note that $m_{G_{N+1}}>m_{G_n}$. From the numerator in \eqref{Tempeqn}, it follows that 
	\begin{align*}
	\Phi_{G_{N+1}}(v,k_1,k_2)&=4m_{G_{N+1}}^2\tau_{G_{N+1}}\left(\frac{\phi_{G_{N+1}}(v)}{4m_{G_{N+1}}^2\tau_{G_{N+1}}}+\phi_1(k_1,k_2)+\frac{\phi_2(k_1,k_2)}{6m_{G_{N+1}}}+\frac{\phi_3(k_1,k_2)}{6m_{G_{N+1}}^2}\right)\\
	&> 4m_{G_{N+1}}^2\tau_{G_{N+1}}\left(\frac{\phi_{G_{N+1}}(v)}{4m_{G_{N+1}}^2\tau_{G_{N+1}}}+\phi_1(k_1,k_2)+\frac{\phi_2(k_1,k_2)}{6m_{G_{N}}}+\frac{\phi_3(k_1,k_2)}{6m_{G_{N}}^2}\right)\\
	&\geq 4m_{G_{N+1}}^2\tau_{G_{N+1}}\left(\frac{\phi_{G_{N}}(v)}{4m_{G_{N}}^2\tau_{G_{N}}}+\phi_1(k_1,k_2)+\frac{\phi_2(k_1,k_2)}{6m_{G_{N}}}+\frac{\phi_3(k_1,k_2)}{6m_{G_{N}}^2}\right)\\
	&=4m_{G_{N+1}}^2\tau_{G_{N+1}}\frac{\Phi_{G_{N}}(v,k_1,k_2)}{4m_{G_{N}}^2\tau_{G_{N}}}>0.
	\end{align*}
	Note that the first inequality is obtained by \eqref{minimum of phi_2 for k} and \eqref{minimum of phi_3}.
\end{proof}

\begin{example}
	{\rm Let $(k_1,k_2)=(1,2)$. Consider $\cG^v=(P_n)^v$ where for each $n\geq 2$, $v$ is a pendent vertex of $P_n$. Examining Example \ref{Ex:asymp for P_n}, it can be seen that $\Phi_{P_{4}}(v,1,2)<0$ and $\Phi_{P_{5}}(v,1,2)>0$; and $\frac{\phi_{P_n}(v)}{4m_{P_n}^2\tau_{P_n}}$ is strictly increasing for $n\geq 2$. Therefore, by Proposition \ref{Proposition:monotone}, $P_n$ is $(v,1,2)$-paradoxical for $n\geq 5$. 
	}
\end{example}

Here is the main result in this section.

\begin{theorem}\label{Theorem:asymtotically paradoxical}
	Given a sequence $\cG^v=(G_n)^v$, $\cG^v$ is asymptotically paradoxical if and only if $\frac{\phi_{G_n}(v)}{4m_{G_n}^2\tau_{G_n}}\rightarrow\infty$ as $n\rightarrow\infty$.
\end{theorem}
\begin{proof}
	We shall prove the sufficiency by contrapositive. Suppose that $\frac{\phi_{G_n}(v)}{4m_{G_n}^2\tau_{G_n}}$ is bounded, say $0<\frac{\phi_{G_n}(v)}{4m_{G_n}^2\tau_{G_n}}\leq L$ for any $n\geq 2$ and for some $L>0$. Then, from \eqref{Tempeqn}, we have 
	\begin{align*}
	&\kappa(\widehat{G}_n(v,k_1,k_2))-\kappa(\widetilde{G}_n(v,k_1,k_2))\\
	<&\frac{\phi_{G_n}(v)+4m_{G_n}^2\tau_{G_n}\phi_1(k_1,k_2)+\frac{2m_{G_n}\tau_{G_n}}{3}\phi_2(k_1,k_2)+\frac{2\tau_{G_n}}{3}\phi_3(k_1,k_2)}{4m_{G_n}^2\tau_{G_n}}\\
	\leq&L+\phi_1(k_1,k_2)+\frac{\phi_2(k_1,k_2)}{6m_{G_n}}+\frac{\phi_3(k_1,k_2)}{6m_{G_n}^2}.
	\end{align*}
	Considering Remark \ref{Remark:phi_1 decrease}, \eqref{minimum of phi_2 for k}, and \eqref{minimum of phi_3}, there exist integers $K_1\geq 0$ and $K_2\geq 0$ with $K_1+K_2\geq 2$ such that $\kappa(\widehat{G}_n(v,K_1,K_2))-\kappa(\widetilde{G}_n(v,K_1,K_2))<0$ for all $n\geq 2$---that is, $\cG^v$ is not asymptotically paradoxical.

	Suppose that $\frac{\phi_{G_n}(v)}{4m_{G_n}^2\tau_{G_n}}$ diverges to infinity. Fix $k_1,k_2\geq 0$ with $k_1+k_2\geq 2$. Since $G_n$ is connected for all $n\geq 1$, $m_{G_n}$ goes to infinity as $n\rightarrow\infty$. It follows from \eqref{Tempeqn} that
	\begin{align*}
	\lim_{n\rightarrow\infty}\left(\kappa(\widehat{G}_n(v,k_1,k_2))-\kappa(\widetilde{G}_n(v,k_1,k_2))\right)=\infty. 
	\end{align*}
	Therefore, $\cG^v$ is asymptotically paradoxical.
\end{proof}

\begin{example}\label{Example:asymp. paradoxical}
	{\rm Revisit Examples \ref{Ex:asymp for K_n}--\ref{Ex:asymp for star}. One can verify that as $n\rightarrow\infty$, we have $\frac{\phi_{K_n}(v)}{4m_{K_n}^2\tau_{K_n}}\rightarrow0$; $\frac{\phi_{C_n}(v)}{4m_{C_n}^2\tau_{C_n}}\rightarrow\infty$; $\frac{\phi_{P_n}(v)}{4m_{P_n}^2\tau_{P_n}}\rightarrow\infty$; $\frac{\phi_{S_n}(v)}{4m_{S_n}^2\tau_{S_n}}\rightarrow 2$ where $v$ is a pendent vertex of $S_n$; and $\frac{\phi_{S_n}(v)}{4m_{S_n}^2\tau_{S_n}}\rightarrow 0$ where $v$ is the centre vertex of $S_n$. By Theorem \ref{Theorem:asymtotically paradoxical}, the sequences $(C_n)^v$ and $(P_n)^v$ are asymptotically paradoxical.}
\end{example}

Now we shall construct $(v,k_1,k_2)$-paradoxical graphs from a connected graph that is not $(v,k_1,k_2)$-paradoxical, by using an asymptotically paradoxical sequence. Given a connected graph $G$ with a vertex $v$, suppose that $G$ is not $(v,k_1,k_2)$-paradoxical. Adding new vertices and edges to $G$, we shall make the resulting graph $(v,k_1,k_2)$-paradoxical.

Note that the case of the equality in the following proposition is used in Section \ref{Section:Asymptotic for trees}.

\begin{proposition}\label{Proposition:dRd with cut vertex}
	Let $G$ be a connected graph, and $v$ be a cut-vertex. Suppose that there are $\ell$ branches $B_1,\dots,B_\ell$ of $G$ at $v$. Then,
	\begin{align}\label{drd:branches}
	\mathbf{d}_G^TQ_{G,v}\mathbf{d}_G=\sum_{k=1}^{\ell}\tau_{G'_k}\mathbf{d}_{B_k}^TQ_{B_k,v}\mathbf{d}_{B_k}\geq \sum_{k=1}^{\ell}\mathbf{d}_{B_k}^TQ_{B_k,v}\mathbf{d}_{B_k}
	\end{align}
	where $G'_k$ is the subgraph induced by $V(G)\backslash (V(B_k)-\{v\})$. This implies that $\phi_G(v)\geq\sum_{k=1}^{\ell}\phi_{B_k}(v)$. Moreover, the two sides are equal if and only if $G$ is a tree.
\end{proposition}
\begin{proof}
	Let $Q_{G,v}=[q_{i,j}]$. By Remark \ref{Remark:entries r}, if $i=v$ or $j=v$, then $q_{i,j}=0$. Consider $i\neq v$ and $j\neq v$. Suppose that $i\in V(B_{k_1})$ and $j\in V(B_{k_2})$ for $k_1\neq k_2$. Since $v$ is a cut-vertex of $G$, we find from Remark \ref{Remark:entries r} that $q_{i,j}=0$. Hence, for $k=1,\dots,\ell$, the submatrix of $Q_{G,v}$ whose rows and columns are indexed by $V(B_k)$ and $V(G)\backslash V(B_k)$, respectively, is the zero matrix. For $k=1,\dots,\ell$, assume $i,j\in V(B_k)$. Since $v$ is a cut-vertex, by Remark \ref{Remark:entries r}, we have $|\mathcal{F}_{G}(i,j;v)|=\tau_{G'_k}|\mathcal{F}_{B_k}(i,j;v)|$ where $G'_k$ is the subgraph induced by $V(G)\backslash (V(B_k)-\{v\})$, with equality if and only if $G'_k$ is a tree. Therefore, the submatrix of $Q_{G,v}$ whose rows and columns are indexed by the vertex set $V(B_k)$ is $\tau_{G'_k}Q_{B_k,v}$.
	
	Let $1\leq k \leq \ell$. For $\mathbf{d}_{B_k}=(d_i)_{i\in V(B_k)}$, let $\widehat{\mathbf{d}}_{B_k}=(\hat{d}_i)_{i\in V(G)}$ where $\hat{d}_i=d_i$ if $i\in V(B_k)$, and $\hat{d}_i=0$ if $i\in V(G)\backslash V(B_k)$. Then, for $1\leq k_1,k_2\leq \ell$,
	\begin{align*}
	\widehat{\mathbf{d}}_{B_{k_1}}^TQ_{G,v}\widehat{\mathbf{d}}_{B_{k_2}}=\mathbf{d}_{B_{k_1}}^T\widetilde{Q}_{G,v}\mathbf{d}_{B_{k_2}}
	\end{align*}
	where $\widetilde{Q}_{G,v}$ is the submatrix of $Q_{G,v}$ whose rows and columns are indexed by $V(B_{k_1})$ and $V(B_{k_2})$, respectively. If $k_1\neq k_2$ then $\mathbf{d}_{B_{k_1}}^T\widetilde{Q}_{G,v}\mathbf{d}_{B_{k_2}}=0$. Furthermore, $\mathbf{d}_{B_{k_1}}^T\widetilde{Q}_{G,v}\mathbf{d}_{B_{k_1}}=\tau_{G'_{k_1}}\mathbf{d}_{B_{k_1}}^TQ_{B_{k_1},v}\mathbf{d}_{B_{k_1}}$. Evidently, $\mathbf{d}_G=\sum_{k=1}^{\ell}\widehat{\mathbf{d}}_{B_k}$. Therefore, the desired result follows.
\end{proof}

\begin{proposition}\label{Proposition:contruction of paradoxical graphs}
	Let $H_i$ be a connected graph with a vertex $v_i$ for $i=1,\dots,\ell$. Suppose that a sequence $\cG^v=(G_n)^v$ is asymptotically paradoxical. Consider a sequence $(\cG')^v=(G'_n)^v$ where for $1\leq n \leq \sum_{i=1}^{\ell}|V(H_i)|$, $G'_n=G_n$, and for $n>\sum_{i=1}^{\ell}|V(H_i)|$, $G'_n$ is the graph obtained from $H_1,\dots,H_\ell$ and $G_{n-\sum_{i=1}^{\ell}|V(H_i)|}$ by identifying the vertices $v_1,\dots,v_\ell,v$. Then, $(\cG')^v$ is asymptotically paradoxical.
\end{proposition}
\begin{proof}
	Suppose that $n>\sum_{i=1}^{\ell}|V(H_i)|$. Let $n_0=n-\sum_{i=1}^{\ell}|V(H_i)|$. Since $v$ is a cut-vertex in $G'_n$, we have $\tau_{G'_n}=\tau_{G_{n_0}}\tau_{H_1}\cdots\tau_{H_\ell}$. Using Proposition \ref{Proposition:dRd with cut vertex}, we obtain
	\begin{align*}
	\frac{\phi_{G'_n}(v)}{4m_{G'_n}^2\tau_{G'_n}}\geq\frac{\phi_{G_{n_0}}(v)+\sum_{i=1}^{\ell}\phi_{H_i}(v_i)}{4(m_{G_{n_0}}+\sum_{i=1}^{\ell}m_{H_i})^2\tau_{G_{n_0}}\tau_{H_1}\cdots\tau_{H_\ell}}.
	\end{align*}
	As $n\rightarrow \infty$, we have $n_0\rightarrow\infty$. Since $(\cG)^v$ is asymptotically paradoxical, by Theorem \ref{Theorem:asymtotically paradoxical} we obtain $\frac{\phi_{G_{n_0}}(v)}{4m_{G_{n_0}}^2\tau_{G_{n_0}}}\rightarrow\infty$ as $n\rightarrow\infty$. It follows that $\frac{\phi_{G'_n}(v)}{4m_{G'_n}^2\tau_{G'_n}}\rightarrow\infty$ as $n\rightarrow\infty$. Therefore, $(\cG')^v$ is asymptotically paradoxical.
\end{proof}
\begin{remark}\label{Remark:discussion for construction}
	{\rm Let a sequence $\cG^v=(G_n)^v$ be asymptotically paradoxical. Suppose that a connected graph $H$ with a vertex $w$ is not $(w,k_1,k_2)$-paradoxical for some integers $k_1$ and $k_2$ with $k_1+k_2\geq 2$. Proposition \ref{Proposition:contruction of paradoxical graphs} tells that regardless of the number of branches of $H$ at $w$, we can obtain a $(v,k_1,k_2)$-paradoxical graph from $H$ by identifying $w$ and the vertex $v$ of $G_n$ for sufficiently large order $n$.}
\end{remark}


\begin{example}
	{\rm Adopting the notation in Remark \ref{Remark:discussion for construction}, consider the following graph $H$:
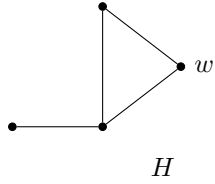
		\begin{center}
			\begin{tikzpicture}[scale=0.80]
			\tikzset{enclosed/.style={draw, circle, inner sep=0pt, minimum size=.10cm, fill=black}}
			
			\node[enclosed, label={left, yshift=0cm: }] (v_1) at (-1,1) {};
			\node[enclosed, label={left, yshift=0cm: }] (v_2) at (-1,-1) {};
			\node[enclosed, label={right, yshift=0cm: $w$}] (v_3) at (0.3,0) {};
			\node[enclosed, label={right, yshift=0cm: }] (v_4) at (-2.5,-1) {};
			\node[label={below, yshift=0cm: $H$}] (G) at (0,-1.2) {};
			\draw (v_1) -- (v_2);
			\draw (v_2) -- (v_3);
			\draw (v_1) -- (v_3);
			\draw (v_2) -- (v_4);
			\end{tikzpicture}
		\end{center}
		One can check from computation that $\phi_{H}(w)=2\mathbf{d}_{H}^TQ_{H,w}\mathbf{d}_{H}=118$ and $\Phi_{H}(w,1,2)<0$. So, $H$ is not $(w,1,2)$-paradoxical. From Example \ref{Example:asymp. paradoxical}, $\cG^v=(P_n)^v$ is asymptotically paradoxical. For ease of exposition, we assume that for each $n\geq 2$, $v$ is a pendent vertex of $P_n$. Suppose that $G'_{n}$ is the graph obtained from $H$ and $P_n$ by identifying $w$ and $v$ as $v$. As discussed in Remark \ref{Remark:discussion for construction}, there must be some $N_0>0$ such that $G'_n$ is $(v,1,2)$-paradoxical for all $n\geq N_0$. We shall minimize such an $N_0$. Using Proposition \ref{Proposition:dRd with cut vertex} and Example \ref{Example:asymp. paradoxical}, we have 
		$$
		\phi_{G'_n}(v)=2\mathbf{d}_{G'_n}^TQ_{G'_n,v}\mathbf{d}_{G'_n}=2\mathbf{d}_{H}^TQ_{H,w}\mathbf{d}_{H}+6\mathbf{d}_{P_n}^TQ_{P_n,v}\mathbf{d}_{P_n}=118+2(n-1)(2n-1)(2n-3).
		$$
		By computation of $\Phi_{G'_n}(v,1,2)$ for $n=2,\dots,5$, $G'_2$, $G'_3$, and $G'_4$ are not $(v,1,2)$-paradoxical, and $G'_5$ is $(v,1,2)$-paradoxical. Furthermore, it can be checked that $\frac{\phi_{G'_n}(v)}{4m_{G'_n}^2\tau_{G'_n}}$ is strictly increasing for $n\geq 5$. Hence, by Proposition \ref{Proposition:monotone}, $G'_n$ is $(v,1,2)$-paradoxical for all $n\geq 5$. In other words, we can construct a $(v,1,2)$-paradoxical graph from $H$ and a path of length at least $4$ by identifying $w$ and a pendent vertex of the path.}  
\end{example}

\section{Asymptotically paradoxical sequences of trees}\label{Section:Asymptotic for trees}
We begin with presenting an outline of this section. Throughout this section, we shall consider sequences $\cG^v=(\cT_n)^v$ of trees, where for each $n\geq 2$, $\cT_n$ is obtained from $\cT_{n-1}$ by an addition of a new pendent vertex or a subdivision of an edge. Then, we examine asymptotic behaviour of such trees upon the addition of twin pendent paths; specifically, we investigate under what circumstances the sequences are asymptotically paradoxical. Considering Theorem \ref{Theorem:asymtotically paradoxical}, we need to understand $\frac{\phi_{\cT_n}(v)}{4m_{\cT_n}^2\tau_{\cT_n}}$. Recall that $\phi_{\cT_n}(v)=2\mathbf{d}_{\cT_n}^TQ_{\cT_n,v}\mathbf{d}_{\cT_n}$. To consider conditions for $\frac{\phi_{\cT_n}(v)}{4m_{\cT_n}^2\tau_{\cT_n}}$ to diverge to infinity, we shall find the minimum of $\mathbf{d}_{\cT_n}^TQ_{\cT_n,v}\mathbf{d}_{\cT_n}$, provided the number of branches of $\cT_n$ at $v$ and the eccentricity of $v$ in each branch are given (Proposition \ref{Proposition:lowest bound of drd for tree}). With the minimum, we provide some condition in terms of the eccentricity of $v$ in $\cT_n$ and the number of branches of $T_n$ at $v$ satisfying some property, in order for the sequences to be asymptotically paradoxical (Theorem \ref{Theorem: asymp. l_n and e(v) growing}).

Here is a sketch of two steps to find the minimum of $\mathbf{d}_{\cT_n}^TQ_{\cT_n,v}\mathbf{d}_{\cT_n}$. By Proposition \ref{Proposition:dRd with cut vertex}, we only need to understand the minimum of $\mathbf{d}_{B}^TQ_{B,v}\mathbf{d}_{B}$ where $B$ is a branch of $\cT_n$ at $v$---that is, the minimum of $\mathbf{d}_{\cT}^TQ_{\cT,v}\mathbf{d}_{\cT}$ where $\cT$ is a tree with a pendent vertex $v$ and the eccentricity of $v$ is given. This minimum is provided in \eqref{lowest bound of dRd for tree} at the end of Step 1. By Proposition \ref{Proposition:dRd with cut vertex}, we establish our desired result in Proposition \ref{Proposition:lowest bound of drd for tree} in Step 2.
\subsection*{Step 1}\label{Subsec:step1}
\begin{figure}[h!]
	\begin{center}
		\begin{tikzpicture}
		\tikzset{enclosed/.style={draw, circle, inner sep=0pt, minimum size=.10cm, fill=black}}
		
		\node[enclosed, label={below, yshift=0cm: $v$}] (v_1) at (0,1) {};
		\node[enclosed, label={below, yshift=0cm: $v_1$}] (v_2) at (1.5,1) {};
		\node[enclosed, label={below, yshift=0cm: $v_2$}] (v_3) at (3,1) {};
		\node[enclosed, label={below, yshift=0cm: $v_3$}] (v_4) at (4.5,1) {};
		\node[enclosed, label={below, yshift=0cm: $v_4$}] (v_5) at (6,1) {};
		\node[enclosed] (v_6) at (1,1.5) {};
		\node[enclosed] (v_10) at (1.5,1.5) {};
		\node[enclosed] (v_7) at (2,1.5) {};
		\node[enclosed] (v_8) at (3,1.5) {};
		\node[enclosed] (v_9) at (3,2) {};
		
		\draw (v_1) -- (v_2);
		\draw (v_2) -- (v_3);
		\draw (v_3) -- (v_4);
		\draw (v_4) -- (v_5);
		\draw (v_2) -- (v_6);
		\draw (v_2) -- (v_7);
		\draw (v_2) -- (v_10);
		\draw (v_3) -- (v_8);
		\draw (v_8) -- (v_9);
		\end{tikzpicture}
	\end{center}
	\caption{}\label{Figure:Example splitting tree}
\end{figure}

Let $\mathcal{T}$ be a tree of order $n$, and $v$ be a pendent vertex in $\mathcal{T}$. Suppose that $\alpha$ is the eccentricity $e_{\cT}(v)$ of $v$ in $\cT$. Then, there exists the path $P=(v_0,v_1,\dots,v_\alpha)$ of length $\alpha$ in $\cT$ where $v_0=v$. Evidently, $v_0$ and $v_\alpha$ are pendent vertices in $\cT$. Let $\mathcal{T}_0$ and $\mathcal{T}_\alpha$ be the trees where $V(\mathcal{T}_0)=\{v_0\}$ and $V(\mathcal{T}_\alpha)=\{v_\alpha\}$. For $k=1,\dots,\alpha-1$, if there are more than two branches of $\cT$ at $v_k$, then we define $\mathcal{T}_k$ to be the tree obtained from $\mathcal{T}$ by deleting two branches except $v_k$ where one contains $v_{k-1}$ and the other $v_{k+1}$; if there are exactly two branches of $\cT$ at $v_k$, then we define $\cT_k$ to be the tree with $V(\cT_k)=\{v_k\}$. Then, $V(\cT_0),\dots,V(\cT_\alpha)$ are mutually disjoint sets. Moreover, for each $k=0,\dots,\alpha$, we have $e_{\mathcal{T}_k}(v_k)\leq \alpha-k$. As an example, if $\cT$ is the tree in Figure \ref{Figure:Example splitting tree}, then $V(\cT_0)=\{v\}$, $V(\cT_3)=\{v_3\}$, and $V(\cT_4)=\{v_4\}$; furthermore, $\cT_1$ and $\cT_2$ are $S_4$ and $P_3$, respectively.

Let $Q_{\mathcal{T},v}=[q_{i,j}]$. Recall that $q_{i,j}=|\mathcal{F}_{\cT}(i,j;v)|$ is the number of $2$-tree spanning forests of $\cT$ having $i$, $j$ in one tree and $v$ in the other. Note that $v=v_0$. In order to understand the structure of $Q_{\mathcal{T},v}$, we shall consider two cases: \begin{enumerate*}[label=(\roman*)]
	\item $i$ and $j$ are in different subtrees; and
	\item $i$ and $j$ are in the same subtree.
\end{enumerate*} Suppose that $i\in V(\mathcal{T}_{k_1})$ and $j\in V(\mathcal{T}_{k_2})$ where $0\leq k_1< k_2\leq \alpha$. For each forest in $\mathcal{F}_{\cT}(i,j;v)$, since  $i$ and $j$ belong to the same subtree in the forest, the subtree must contain $v_{k_1}$ and $v_{k_2}$. For any vertex $w$ on the subpath of $\cT$ with $i$ and $j$ as the pendent vertices, we have $\mathrm{dist}_{\mathcal{T}}(v,v_{k_1})\leq \mathrm{dist}_{\mathcal{T}}(v,w)$. Hence, by Remark \ref{Remark:entries r}, $q_{i,j}=k_1$ for $i\in V(\mathcal{T}_{k_1})$ and $j\in V(\mathcal{T}_{k_2})$ with $0\leq k_1<k_2\leq \alpha$. 

Assume that $i,j$ are in $V(\cT_k)$ for some $1\leq k\leq \alpha$. Consider the subpath $P'$ of $\cT_k$ with $i$ and $j$ as the pendent vertices. Suppose that $w_0$ is the vertex on $P'$ such that $\mathrm{dist}_{\mathcal{T}_k}(v_k,w_0)\leq \mathrm{dist}_{\mathcal{T}_k}(v_k,w)$ for $w\in V(P')$. Then, $\mathrm{dist}_{\mathcal{T}}(v,w_0)=k+\mathrm{dist}_{\mathcal{T}_k}(v_k,w_0)$. Let  $Q_{\mathcal{T}_k,v_k}=[\tilde{q}_{i,j}]$. By Remark \ref{Remark:entries r}, we have $q_{i,j}=k+\tilde{q}_{i,j}$.

Labelling the rows and columns of $Q_{\mathcal{T},v}$ in order of $v,V(\mathcal{T}_1),\dots,V(\mathcal{T}_\alpha)$, we obtain the following structure:
$$
Q_{\mathcal{T},v}=
\left[\begin{array}{c|c|c|c|c|c}
0 & 0 & 0 & 0 & \cdots & 0\\\hline
0 & J+Q_{\mathcal{T}_1,v_1} & J & J & \cdots & J\\\hline
0 & J & 2J+Q_{\mathcal{T}_2,v_2} & 2J & \cdots & 2J\\\hline
0 & J & 2J & 3J+Q_{\mathcal{T}_3,v_3} & \cdots & \vdots\\\hline
\vdots & \vdots & \vdots & \vdots & \ddots & (\alpha-1)J\\\hline
0 & J & 2J & 3J & \cdots & \alpha J+Q_{\mathcal{T}_\alpha,v_\alpha}
\end{array}\right]
$$
where the $J$s in the blocks of $Q_{\cT,v}$ are appropriately sized. Let $n_k=|V(\mathcal{T}_k)|$ for $k=0,\dots,\alpha$. Note that $n_0=n_\alpha=1$. Then, $Q_{\mathcal{T},v}$ can be recast as 
\begin{align*}
Q_{\mathcal{T},v}&=\sum_{i=0}^{\alpha-1}\begin{bmatrix}
0 & 0 \\
0 & J_{n-(n_0+\cdots+n_i)}
\end{bmatrix}+\mathrm{diag}(0,Q_{\mathcal{T}_1,v_1},\dots,Q_{\mathcal{T}_\alpha,v_\alpha})\\
&=\sum_{i=0}^{\alpha-1}\begin{bmatrix}
\mathbf{0}_{n_0+\cdots+n_i}\\
\mathbf{1}_{n-(n_0+\cdots+n_i)}
\end{bmatrix}\begin{bmatrix}
\mathbf{0}_{n_0+\cdots+n_i}^T & \mathbf{1}_{n-(n_0+\cdots+n_i)}^T
\end{bmatrix}+\mathrm{diag}(0,Q_{\mathcal{T}_1,v_1},\dots,Q_{\mathcal{T}_\alpha,v_\alpha})
\end{align*}
where $n=n_0+n_1+\dots+n_\alpha$. 

Now, we shall compute $\mathbf{d}^T_\mathcal{T}Q_{\mathcal{T},v}\mathbf{d}_\mathcal{T}$. Let $\bx^T=\begin{bmatrix}
0 & \mathbf{d}_{\mathcal{T}_1}^T & \cdots & \mathbf{d}_{\mathcal{T}_{\alpha-1}}^T & 0
\end{bmatrix}$ and $\by=\mathbf{e}_{v}+\sum_{i=1}^{\alpha-1}2\mathbf{e}_{v_i}+\mathbf{e}_{v_\alpha}$. Then,
\begin{align*}
\bx^T Q_{\mathcal{T},v}\bx&=\sum_{i=0}^{\alpha-2}\left(\mathbf{d}_{\cT_{i+1}}^T\mathbf{1}+\cdots +\mathbf{d}_{\cT_{\alpha-1}}^T\mathbf{1}\right)^2+\sum_{i=1}^{\alpha-1}\mathbf{d}^T_{\mathcal{T}_i}Q_{\mathcal{T}_i,v_i}\mathbf{d}_{\mathcal{T}_i}\\
&=4\sum_{i=0}^{\alpha-2}\left(\sum_{j=i+1}^{\alpha-1}(n_j-1)\right)^2+\sum_{i=1}^{\alpha-1}\mathbf{d}^T_{\mathcal{T}_i}Q_{\mathcal{T}_i,v_i}\mathbf{d}_{\mathcal{T}_i}=4\sum_{i=1}^{\alpha-1}\left(\sum_{j=i}^{\alpha-1}(n_j-1)\right)^2+\sum_{i=1}^{\alpha-1}\mathbf{d}^T_{\mathcal{T}_i}Q_{\mathcal{T}_i,v_i}\mathbf{d}_{\mathcal{T}_i}.
\end{align*}

We can find that the submatrix of $Q_{\cT,v}$ whose rows and columns are indexed by $\{v_0,\dots,v_\alpha\}$ is $[\mathrm{min}(i,j)]_{0\leq i,j\leq \alpha}$. So, $\left(\sum_{k=0}^\alpha\be_{v_k}\right)^TQ_{\cT,v}\left(\sum_{k=0}^\alpha\be_{v_k}\right)$ is the sum of all entries in $[\mathrm{min}(i,j)]_{0\leq i,j\leq \alpha}$. Thus, from $\by=2\left(\sum_{k=0}^\alpha\be_{v_k}\right)-(\be_v+\be_\alpha)$, we have
\begin{align*}
\by^TQ_{\mathcal{T},v}\by&=4\mathbf{1}^T[\mathrm{min}(i,j)]_{0\leq i,j\leq \alpha}\mathbf{1}-4(\be_v+\be_\alpha)^TQ_{\mathcal{T},v}\left(\sum_{k=0}^\alpha\be_{v_k}\right)+(\be_v+\be_\alpha)^TQ_{\mathcal{T},v}(\be_v+\be_\alpha)\\
&=\frac{2}{3}\alpha(\alpha+1)(2\alpha+1)-2\alpha(\alpha+1)+\alpha=\frac{1}{3}\alpha(2\alpha-1)(2\alpha+1).
\end{align*}

Finally, we find
\begin{align*}
&\sum_{i=0}^{\alpha-1}\bx^T\begin{bmatrix}
0 & 0 \\
0 & J_{n-(n_0+\cdots+n_i)}
\end{bmatrix}\by\\
=&\sum_{i=0}^{\alpha-2}\left(\mathbf{d}_{\cT_{i+1}}^T\mathbf{1}+\cdots +\mathbf{d}_{\cT_{\alpha-1}}^T\mathbf{1}\right)\left(2(\alpha-i)-1\right)\\
=&2\sum_{i=0}^{\alpha-2}\left(\sum_{j=i+1}^{\alpha-1}(n_j-1)\right)\left(2(\alpha-i)-1\right)=2\sum_{i=1}^{\alpha-1}\left(\sum_{j=i}^{\alpha-1}(n_j-1)\right)\left(2(\alpha-i)+1\right).
\end{align*}
From Remark \ref{Remark:entries r}, for each $k=0,\dots,\alpha$, we have $|\mathcal{F}_{\cT_k}(l,v_k;v_k)|=0$ for $l\in V(\cT_k)$. So, the $v_k^\text{th}$ column of $\mathrm{diag}(0,Q_{\mathcal{T}_1,v_1},\dots,Q_{\mathcal{T}_\alpha,v_\alpha})$ is the zero vector. This implies $\bx^T\mathrm{diag}(0,Q_{\mathcal{T}_1,v_1},\dots,Q_{\mathcal{T}_\alpha,v_\alpha})\by=0$. Hence, 
\begin{align*}
2\bx^TQ_{\mathcal{T},v}\by&=4\sum_{i=1}^{\alpha-1}\left(\sum_{j=i}^{\alpha-1}(n_j-1)\right)(2(\alpha-i)+1).
\end{align*}
Note that $\bd_\cT=\bx+\by$. Therefore, for a tree $\cT$ with a pendent vertex $v$,
\begin{align}\label{identity:dRd for tree with v}
\begin{split}
\mathbf{d}^T_\mathcal{T}Q_{\mathcal{T},v}\mathbf{d}_\mathcal{T}=&\bx^T Q_{\mathcal{T},v}\bx+2\bx^TQ_{\mathcal{T},v}\by+\by^TQ_{\mathcal{T},v}\by\\
=&4\sum_{i=1}^{\alpha-1}\left(\sum_{j=i}^{\alpha-1}(n_j-1)\right)^2+\sum_{i=1}^{\alpha-1}\mathbf{d}^T_{\mathcal{T}_i}Q_{\mathcal{T}_i,v_i}\mathbf{d}_{\mathcal{T}_i}\\
&+4\sum_{i=1}^{\alpha-1}\left(\sum_{j=i}^{\alpha-1}(n_j-1)\right)(2(\alpha-i)+1)+\frac{1}{3}\alpha(2\alpha-1)(2\alpha+1).
\end{split}
\end{align}
\begin{figure}[h!]
	\begin{center}
		\begin{tikzpicture}
		\tikzset{enclosed/.style={draw, circle, inner sep=0pt, minimum size=.10cm, fill=black}}
		
		\node[enclosed, label={below, yshift=0cm: $v$}] (v_1) at (0,1) {};
		\node[enclosed, label={below, yshift=0cm: $v_1$}] (v_2) at (1,1) {};
		\node[enclosed, label={below, yshift=0cm: $v_2$}] (v_3) at (2,1) {};
		\node[enclosed, label={below, yshift=0cm: $v_{\alpha-1}$}] (v_4) at (3,1) {};
		\node[enclosed, label={below, yshift=0cm: $v_\alpha$}] (v_5) at (4,1) {};
		\node[enclosed] (v_6) at (2.1340,1.5) {};
		\node[enclosed] (v_7) at (2.5,1.8660) {};
		\node[enclosed] (v_8) at (3.5,1.8660) {};
		\node[enclosed] (v_9) at (3.8660,1.5) {};
		\node[label={center, yshift=0cm: $\dots$}] (v_10) at (3,2) {};
		
		\draw (v_1) -- (v_2);
		\draw (v_2) -- (v_3);
		\draw[thick, loosely dotted] (v_3) -- (v_4);
		\draw (v_4) -- (v_5);
		\draw (v_4) -- (v_6);
		\draw (v_4) -- (v_7);
		\draw (v_4) -- (v_8);
		\draw (v_4) -- (v_9);
		\end{tikzpicture}
	\end{center}
	\caption{A broom on $n$ vertices with exactly $(n-\alpha)$ pendent vertices having a common neighbour.}\label{Figure:Broom}
\end{figure}
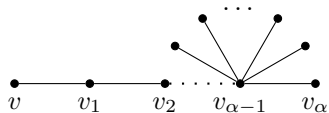

\begin{example}\label{Ex:Broom dRd}
	{\rm Let $n\geq\alpha\geq 1$, and $\mathcal{B}_{n,\alpha}$ be the broom with vertices $v,v_1,\dots,v_\alpha$ in Figure \ref{Figure:Broom}. Let $v_0=v$, and $X=\{0,\dots,\alpha\}\backslash\{\alpha-1\}$. Suppose that for $i\in X$, $\cT_i$ is the tree with $V(\cT_i)=\{v_i\}$, and $\cT_{\alpha-1}$ is the subtree induced by $V(\mathcal{B}_{n,\alpha})\backslash \{v_1,\dots,v_{\alpha-2},v_\alpha\}$. Then, $\cT_{\alpha-1}$ is a star of order $n-\alpha$ with the centre vertex $v_{\alpha-1}$. Let $n_i=|V(\cT_i)|$ for $i=0,\dots,\alpha$. By \eqref{identity:dRd for tree with v} and Example \ref{Ex:asymp for star}, we obtain
		\begin{align*}
		\mathbf{d}_{\mathcal{B}_{n,\alpha},v}^TQ_{\mathcal{B}_{n,\alpha},v}\mathbf{d}_{\mathcal{B}_{n,\alpha},v}=&4\sum_{i=1}^{\alpha-1}(n-\alpha-1)^2+\mathbf{d}^T_{S_{n-\alpha}}Q_{S_{n-\alpha},v_{\alpha-1}}\mathbf{d}_{S_{n-\alpha}}\\
		&+4\sum_{i=1}^{\alpha-1}(n-\alpha-1)(2(\alpha-i)+1)+\frac{1}{3}\alpha(2\alpha-1)(2\alpha+1)\\
		=&4(\alpha-1)(n-\alpha-1)^2+(n-\alpha-1)(4\alpha^2-3)+\frac{1}{3}\alpha(2\alpha-1)(2\alpha+1).
		\end{align*}} 
\end{example}

We continue \eqref{identity:dRd for tree with v} with the same hypotheses and notation. We consider $\mathbf{d}^T_{\mathcal{T}_i}Q_{\mathcal{T}_i,v_i}\mathbf{d}_{\mathcal{T}_i}$ for $i=1,\dots,\alpha-1$ in \eqref{identity:dRd for tree with v}. Note that for $i=1,\dots,\alpha-1$, $v_i$ is not necessarily a pendent vertex in $\cT_{v_i}$. The result for the minimum of $\mathbf{d}^T_{\mathcal{T}_i}Q_{\mathcal{T}_i,v_i}\mathbf{d}_{\mathcal{T}_i}$ appears in the paper \cite{KirklandZeng} as the minimum of $\mathbf{d}^T_{\mathcal{T}_i}(2\mathbf{f}_{\cT_i}^v\mathbf{1}^T-F_{\cT_i})\mathbf{d}_{\mathcal{T}_i}$. We shall introduce the result, which is proved by induction in \cite{KirklandZeng}, with a different proof by using the combinatorial interpretation for entries in $Q_{\mathcal{T}_i,v_i}$.

\begin{lemma}\cite{KirklandZeng}\label{Lemma:lower bound of dRd for tree}
	Let $\cT$ be a tree of order $n\geq 2$ with a vertex $v$. Then,
	\begin{align*}
	\mathbf{d}^T Q_{\cT,v} \mathbf{d}\geq n-1
	\end{align*}
	with equality if and only if for $n=2$, $\cT=P_2$ and for $n\geq 3$, $\cT=S_n$ and $v$ is the centre vertex.
\end{lemma}
\begin{proof}
	Let $Q_{\cT,v}=[q_{i,j}]$. By Remark \ref{Remark:entries r}, we have $q_{ii}=\mathrm{dist}(i,v)\geq 1$ whenever $i\neq v$. The degree of each vertex is at least $1$. So, we have $\mathbf{d}^T Q_{\cT,v} \mathbf{d}\geq (n-1)$. To attain the equality, $q_{i,j}=0$ if $i\neq j$. From Remark \ref{Remark:entries r}, we can find that $v$ is a cut-vertex so that $\cT-v$ consists of $n-1$ isolated vertices. Therefore, our desired result is obtained.
\end{proof}

Applying Lemma \ref{Lemma:lower bound of dRd for tree} to $\mathbf{d}^T_{\mathcal{T}_i}Q_{\mathcal{T}_i,v_i}\mathbf{d}_{\mathcal{T}_i}$ in \eqref{identity:dRd for tree with v} for each $i=1,\dots,\alpha-1$, we obtain $\sum_{i=1}^{\alpha-1}\mathbf{d}^T_{\mathcal{T}_i}Q_{\mathcal{T}_i,v_i}\mathbf{d}_{\mathcal{T}_i}\geq n-\alpha-1$. Thus, $\mathbf{d}^T_\mathcal{T}Q_{\mathcal{T},v}\mathbf{d}_\mathcal{T}$ in \eqref{identity:dRd for tree with v} is bounded below as follows:
\begin{align*}
\begin{split}
\mathbf{d}^T_\mathcal{T}Q_{\mathcal{T},v}\mathbf{d}_\mathcal{T}
\geq& 4\sum_{i=1}^{\alpha-1}\left(\sum_{j=i}^{\alpha-1}(n_j-1)\right)^2+4\sum_{i=1}^{\alpha-1}\left(\sum_{j=i}^{\alpha-1}(n_j-1)\right)(2(\alpha-i)+1)\\
&+(n-\alpha-1)+\frac{1}{3}\alpha(2\alpha-1)(2\alpha+1).
\end{split}
\end{align*}
Consider
\begin{align}\label{temp:iden}
&\sum_{i=1}^{\alpha-1}\left(\sum_{j=i}^{\alpha-1}(n_j-1)\right)^2+\sum_{i=1}^{\alpha-1}\left(\sum_{j=i}^{\alpha-1}(n_j-1)\right)(2(\alpha-i)+1)\\\nonumber
=&\left[(n_1+\cdots+n_{\alpha-1}-(\alpha-1))^2+(n_1+\cdots+n_{\alpha-1}-(\alpha-1))(2\alpha-1)\right]\\\nonumber
&+\left[(n_2+\cdots+n_{\alpha-1}-(\alpha-2))^2+(n_2+\cdots+n_{\alpha-1}-(\alpha-2))(2\alpha-3)\right]\\\nonumber
&+\cdots+\left[(n_{\alpha-1}-1)^2+(n_{\alpha-1}-1)3\right].
\end{align}
Since $n_1+\cdots+n_{\alpha-1}$ is constant, we find that the minimum of \eqref{temp:iden} is attained as $(n-\alpha-1)(n+\alpha-2)$ at $n_1=n-\alpha$ and $n_2=\cdots=n_{\alpha-1}=1$. Therefore, when $v$ is a pendent vertex, we have
\begin{align}\label{lowest bound of dRd for tree}
\mathbf{d}^T_\mathcal{T}Q_{\mathcal{T},v}\mathbf{d}_\mathcal{T}\geq (n-\alpha-1)(4n+4\alpha-7)+\frac{1}{3}\alpha(2\alpha-1)(2\alpha+1)
\end{align}
where equality holds if and only if $\cT$ is a broom $\mathcal{B}_{n,\alpha}$ with $v,v_1,\dots,v_\alpha$ described below:
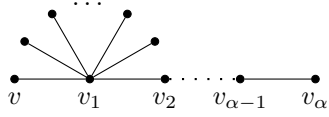
\begin{center}
	\begin{tikzpicture}
	\tikzset{enclosed/.style={draw, circle, inner sep=0pt, minimum size=.10cm, fill=black}}
	
	\node[enclosed, label={below, yshift=0cm: $v$}] (v_1) at (0,1) {};
	\node[enclosed, label={below, yshift=0cm: $v_1$}] (v_2) at (1,1) {};
	\node[enclosed, label={below, yshift=0cm: $v_2$}] (v_3) at (2,1) {};
	\node[enclosed, label={below, yshift=0cm: $v_{\alpha-1}$}] (v_4) at (3,1) {};
	\node[enclosed, label={below, yshift=0cm: $v_{\alpha}$}] (v_5) at (4,1) {};
	\node[enclosed] (v_6) at (0.1340,1.5) {};
	\node[enclosed] (v_7) at (0.5,1.8660) {};
	\node[enclosed] (v_8) at (1.5,1.8660) {};
	\node[enclosed] (v_9) at (1.8660,1.5) {};
	\node[label={center, yshift=0cm: $\dots$}] (v_10) at (1,2) {};
	
	\draw (v_1) -- (v_2);
	\draw (v_2) -- (v_3);
	\draw[thick, loosely dotted] (v_3) -- (v_4);
	\draw (v_4) -- (v_5);
	\draw (v_2) -- (v_6);
	\draw (v_2) -- (v_7);
	\draw (v_2) -- (v_8);
	\draw (v_2) -- (v_9);
	\end{tikzpicture}
\end{center}

\subsection*{Step 2}
The following is the result for the minimum of $\mathbf{d}^T_\mathcal{T}Q_{\mathcal{T},v}\mathbf{d}_\mathcal{T}$ where $\cT$ is a tree with a vertex $v$. (The vertex $v$ is not necessarily a pendent vertex.)

\begin{proposition}\label{Proposition:lowest bound of drd for tree}
	Let $\cT$ be a tree with a vertex $v$. Suppose that $B_1,\dots,B_\ell$ are the branches of $\cT$ at $v$ for some $\ell\geq 1$. Let $n_i=|V(B_i)|$, and let $e_i=e_{B_i}(v)$ for $i=1,\dots,\ell$. Then,
	\begin{align*}
	\mathbf{d}^T_\mathcal{T}Q_{\mathcal{T},v}\mathbf{d}_\mathcal{T}\geq \sum_{i=1}^{\ell}\left[(n_i-e_i-1)(4n_i+4e_i-7)+\frac{1}{3}e_i(2e_i-1)(2e_i+1)\right]
	\end{align*}
	where equality holds if and only if for $i=1,\dots,\ell$, each branch $B_i$ is a broom $\mathcal{B}_{n_i,e_i}$ such that if $n_i>e_i+1$, then $v$ is one of the $(n_i-e_i)$ pendent vertices having a common neighbour; if $n_i=e_i+1$, then $v$ is a pendent vertex in $\mathcal{B}_{n_i,e_i}$ (which is a path).
\end{proposition}
\begin{proof}
	The conclusions can be readily established by Proposition \ref{Proposition:dRd with cut vertex} and \eqref{lowest bound of dRd for tree}.
\end{proof}

Hereafter, the symbols $\omega$, $\mathcal{O}$ and $\Theta$ stand for the small Omega notation, the big O notation and the big Theta notation, respectively (see \cite{Arora:Complexity}). 

As mentioned in the beginning of this section, we consider the following sequence $\cG^v=(\cT_n)^v$ of trees, where $V(\cT_1)=\{v\}$ and for each $n\geq 2$, $\cT_n$ is obtained from $\cT_{n-1}$ by adding a new pendent vertex to $\cT_{n-1}$, or by subdividing an edge in $\cT_{n-1}$ into two edges connecting to a new vertex. We denote by $\alpha_n(x)$ and $\ell_n(x)$ the eccentricity of $x$ in $\cT_n$ and the number of branches of $\cT_n$ at $x$, respectively. For the rest of this section, we use $\alpha_n(\cdot)$ and $\ell_n(\cdot)$ only for the specified vertex $v$ of the trees in the sequence, so we simply write $\alpha_n(v)$ and $\ell_n(v)$ as $\alpha_n$ and $\ell_n$. 

Define $B_1^{(1)}=\cT_1$ and $\ell_1=1$. Assume that for $n\geq 2$, $B_1^{(n-1)},\dots,B_{\ell_{n-1}}^{(n-1)}$ are the branches of $\cT_{n-1}$ at $v$. Let $\{w\}=V(\cT_n)\backslash V(\cT_{n-1})$. Consider the case $\ell_n-\ell_{n-1}=1$. Then, $w$ must be added to the vertex $v$ in $\cT_{n-1}$ to form $\cT_n$. For this case, we define $B_i^{(n)}$ as $B_i^{(n-1)}$ for $i=1,\dots,\ell_n-1$, and define $B_{\ell_n}^{(n)}$ as the path $(v,w)$. Suppose $\ell_n=\ell_{n-1}$. Then, there exists exactly one branch $B_k^{(n-1)}$ for some $k\in\{1,\dots,\ell_{n-1}\}$ such that $w$ is adjacent to at least a vertex of $B_k^{(n-1)}$ in $\cT_n$. We define $B_i^{(n)}$ as $B_i^{(n-1)}$ for $1\leq i\leq \ell_{n-1}$ with $i\neq k$, and define $B_k^{(n)}$ as the induced subtree of $\cT_n$ by $V\left(B_k^{(n-1)}\right)\cup\{w\}$. Hence, we may define $$\beta_{n}=|\{i|e_{B_i^{(k)}}(v)=\Theta(\alpha_k),i=1,\dots,\ell_{n}\}|.$$
Note that $\beta_{n}$ is the number of branches of $\cT_{n}$ at $v$ such that the eccentricity of $v$ in a branch is asymptotically bounded above and below by the eccentricity of $v$ in $\cT_k$.

\begin{example}
	{\rm If $\cG^v=(P_n)^v$ where $v$ is a pendent vertex for $n\geq 2$, then $\alpha_n=n-1$ and $\ell_n=\beta_n=1$. If $\cG^v=(S_n)^v$ where $v$ is the centre vertex for $n\geq 3$, then $\alpha_n=1$ and $\ell_n=\beta_n=n-1$.}
\end{example}

\begin{remark}
	{\rm Consider a sequence $\cG^v=(\cT_n)^v$ of trees. Evidently, $\beta_n\leq \ell_n=\mathcal{O}(n)$ and $\alpha_n=\mathcal{O}(n)$.  Since $\alpha_n=\mathrm{max}\{e_{B_i^{(n)}}(v)|1\leq i\leq \ell_n\}$, we have $\beta_n\geq 1$.}
\end{remark}

Here is the main result in this section.

\begin{theorem}\label{Theorem: asymp. l_n and e(v) growing}
	Let $\cG^v=(\cT_n)^v$ be a sequence of trees. If $\beta_n\alpha_n^3=\omega(n^2)$, then $\cG^v$ is asymptotically paradoxical.
\end{theorem}
\begin{proof}
	Suppose that $\beta_n\alpha_n^3=\omega(n^2)$. For $n\geq 2$, suppose that $B_1^{(n)},\dots,B_{\ell_n}^{(n)}$ are the branches of $\cT_n$ at $v$. Let $e_i^{(n)}=e_{B_i^{(n)}}(v)$ and $k_i^{(n)}=\left|V\left(B_i^{(n)}\right)\right|$ for $i=1,\dots,\ell_n$. We may assume that $e_j^{(n)}=\Theta(\alpha_n)$ for $j=1,\dots,\beta_n$. Then, for each $j=1,\dots,\beta_n$, there exist $C_j>0$ and $N_j>0$ such that $e_j^{(n)}\geq C_j \alpha_n$ for all $n\geq N_j$. Choose $C_0=\mathrm{min}\{C_j|j=1,\dots,\beta_n\}$ and $N_0=\mathrm{max}\{N_j|j=1,\dots,\beta_n\}$. Then, $e_j^{(n)}\geq C_0 \alpha_n$ for all $n\geq N_0$ and $1\leq j\leq \beta_n$. By Proposition \ref{Proposition:lowest bound of drd for tree}, for $n\geq N_0$, we have
	\begin{align*}
	\frac{\phi_{\cT_n}(v)}{4m_{\cT_n}^2\tau_{\cT_n}}&=\frac{2\mathbf{d}^T_{\cT_n}Q_{\mathcal{T}_n,v}\mathbf{d}_{\cT_n}}{4(n-1)^2}\\
	&\geq \frac{\sum_{i=1}^{\ell_n}\left[\left(k_i^{(n)}-e_i^{(n)}-1\right)\left(4k_i^{(n)}+4e_i^{(n)}-7\right)+\frac{1}{3}e_i^{(n)}\left(2e_i^{(n)}-1\right)\left(2e_i^{(n)}+1\right)\right]}{2(n-1)^2}\\
	&\geq \frac{\beta_nC_0\alpha_n(2C_0\alpha_n-1)(2C_0\alpha_n+1)}{6(n-1)^2}.
	\end{align*}
	Since $\beta_n\alpha_n^3=\omega(n^2)$, we have $\frac{\phi_{\cT_n}(v)}{4m_{\cT_n}^2\tau_{\cT_n}}\rightarrow\infty$ as $n$ goes to infinity. Therefore, the conclusion follows.
\end{proof}

\begin{corollary}\label{Corollary:eccentricity and asymp paradoxical}
	Suppose that $\cG^v=(\cT_n)^v$ is a sequence of trees $\cT_n$ such that $\alpha_n=\omega(n^\frac{2}{3})$. Then, $\cG^v$ is asymptotically paradoxical.
\end{corollary}
\begin{proof}
	It is straightforward from Theorem \ref{Theorem: asymp. l_n and e(v) growing}.
\end{proof}

\begin{corollary}\label{Corollary:diameter and asymp paradoxical}
	Suppose that $\cG^v=(\cT_n)^v$ is a sequence of trees $\cT_n$ such that $\mathrm{diam}(\cT_n)=\omega(n^\frac{2}{3})$. Then, $\cG^v$ is asymptotically paradoxical.
\end{corollary}
\begin{proof}
	Let $P$ be a longest path in $\cT_n$. Suppose that $w_0$ is the vertex on $P$ such that $\mathrm{dist}(v,w_0)\leq \mathrm{dist}(v,w)$ for all vertices $w$ on $P$. Then, $\alpha_n\geq \mathrm{dist}(v,w_0)+\frac{1}{2}\mathrm{diam}(\cT_n)$. By Corollary \ref{Corollary:eccentricity and asymp paradoxical}, our desired result follows.
\end{proof}

A \textit{rooted} tree is a tree with a vertex designated as the root such that every edge is directed away from the root. A \textit{leaf} in a rooted tree is a vertex whose degree is $1$. The \textit{depth} of a vertex $v$ in a rooted tree is the distance between $v$ and the root. The \textit{height} of a rooted tree is the maximum distance from the root to all leaves. 

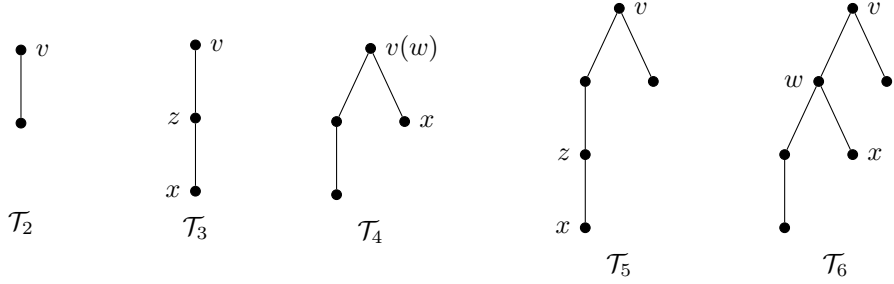
\begin{figure}[h!]
	\hfill
	\adjustbox{valign=c}{
		\begin{forest}
			[$v$
			[]
			]
			\node at (current bounding box.south)
			[below=6.5ex]
			{$\cT_2$};
	\end{forest}}
	\hfill
	\adjustbox{valign=c}{
		\begin{forest}
			[$v$
			[$z$[$x$]]
			]
			\node at (current bounding box.south)
			[below=1ex]
			{$\cT_3$};
	\end{forest}}
	\hfill
	\adjustbox{valign=c}{
		\begin{forest}
			[$v(w)$
			[[]]
			[$x$]
			]
			\node at (current bounding box.south)
			[below=1ex]
			{$\cT_4$};
	\end{forest}}
	\hfill
	\adjustbox{valign=c}{
		\begin{forest}
			[$v$
			[[$z$[$x$]]]
			[]
			]
			\node at (current bounding box.south)
			[below=1ex]
			{$\cT_5$};
	\end{forest}}
	\hfill
	\adjustbox{valign=c}{
		\begin{forest}
			[$v$
			[$w$[[]][$x$]]
			[]
			]
			\node at (current bounding box.south)
			[below=1ex]
			{$\cT_6$};
	\end{forest}}
	\hfill\mbox{}
	\caption{A sequence of rooted trees considered in Example \ref{Example:a sequence of rooted trees}.}\label{Figure:a seq of rooted}
\end{figure}

\begin{example}\label{Example:a sequence of rooted trees}
	{\rm Let $\cG^v=(\cT_n)^v$ be a sequence of trees. For each $n\geq 1$, $\cT_n$ can be considered as a rooted tree at $v$. We may also regard branches $B_1^{(n)},\dots,B_{\ell_n}^{(n)}$ of $\cT_n$ at $v$ as rooted trees at $v$. For each $n\geq 3$, let $\cT_n$ be obtained from $\cT_{n-1}$ as follows: if $e_{B_1^{(n-1)}}(v)=\lfloor n^{c_0}\rfloor-1$, then a new vertex $x$ is added to a leaf $z$ of $B_1^{(n-1)}$ such that the depth of $z$ is the height of $B_1^{(n-1)}$; if $e_{B_1^{(n-1)}}(v)=\lfloor n^{c_0}\rfloor$, then a new vertex $x$ is added to a vertex $w$ in $\cT_{n-1}$ such that $\mathrm{dist}(v,w)<e_{B_1^{(n-1)}}(v)$. Assume that $c_0=0.7$. Considering $\lfloor 3^{c_0}\rfloor=\lfloor 4^{c_0}\rfloor=2$ and $\lfloor 5^{c_0}\rfloor=\lfloor 6^{c_0}\rfloor=3$, one of all possible sequences can be obtained as in Figure \ref{Figure:a seq of rooted}. Note that the very left branch of each rooted tree at $v$ in that figure is $B_1^{(n)}$ for $n=2,\dots,6$. Then, $e_{B_1^{(n)}}(v)\geq e_{B_k^{(n)}}(v)$ for all $n\geq 2$ and $2\leq k\leq\ell_n$. Moreover, $e_{B_1^{(n)}}(v)\geq n^{c_0}-1$ for all $n\geq 2$. By Corollary \ref{Corollary:eccentricity and asymp paradoxical}, $\cG^v$ is asymptotically paradoxical---that is, for integers $k_1,k_2\geq 0$ with $k_1+k_2\geq 2$, $\cT_n$ is $(v,k_1,k_2)$-paradoxical for sufficiently large $n$.}
\end{example}

From the following example, the converses of Theorem \ref{Theorem: asymp. l_n and e(v) growing}, Corollaries \ref{Corollary:eccentricity and asymp paradoxical} and \ref{Corollary:diameter and asymp paradoxical} do not hold.

\begin{example}\label{Example:Counter example broom}
	{\rm Consider a sequence $\cG^v=(\cT_n)^v$ where for $n\geq 4$, $\cT_n$ is a broom $\mathcal{B}_{n,\alpha_n}$ with $\alpha_n\geq 3$. Suppose that for each $n\geq 4$, $v$ is the pendent vertex of $\mathcal{B}_{n,\alpha_n}$ that does not have any common neighbour with other pendent vertices in $\mathcal{B}_{n,\alpha_n}$. Clearly, $\beta_n=1$. Suppose that $\alpha_n=\omega(1)$. By Example \ref{Ex:Broom dRd}, we obtain
		\begin{align*}
		\frac{\phi_{\mathcal{B}_{n,\alpha_n}}(v)}{4m_{\mathcal{B}_{n,\alpha_n}}^2\tau_{\mathcal{B}_{n,\alpha_n}}}=&\frac{4(\alpha_n-1)(n-\alpha_n-1)^2+(n-\alpha_n-1)(4\alpha_n^2-3)+\frac{1}{3}\alpha_n(2\alpha_n-1)(2\alpha_n+1)}{2(n-1)^2}\\
		\geq&\frac{2(\alpha_n-1)(n-\alpha_n-1)^2}{(n-1)^2}
		\end{align*}
		for $n\geq 4$. Since $n^2\alpha_n=\omega(n^2)$, we have $\frac{\phi_{\cT_n}(v)}{4m_{\cT_n}^2\tau_{\cT_n}}\rightarrow\infty$ as $n$ goes to infinity. Therefore, $\cG^v$ is asymptotically paradoxical. Moreover, we have $\beta_n\alpha_n^3=\omega(1)$.}
\end{example}

\bigskip
{\bf Acknowledgment.} The author is grateful to Lorenzo Ciardo at the University of Oslo for discussions that initiated this study, and Steve Kirkland at the University of Manitoba for his encouragement, reviews, and comments.

{\bf Post-acknowledgment.} This article was submitted for publication in the Electronic Journal of Linear Algebra on February 5, 2021. The author later noticed that Proposition \ref{gen:H1andH2} of this manuscript is equivalent to Theorem 2.1 in the article `arXiv:2108.01061' and its publication: N. Faught, M. Kempton, and A. Knudson, A 1-separation formula for the graph Kemeny constant and Braess edges, \textit{Journal of Mathematical Chemistry}, 1--21, 2021.  


\end{document}